\documentclass[11pt]{amsart}
\usepackage{amsmath,amsthm,amsfonts,amscd,amssymb}
\usepackage[mathscr]{eucal}
\usepackage[all]{xy}
\usepackage[colorlinks]{hyperref}

\newtheorem{thm}{Theorem}
\newtheorem{lemma}[thm]{Lemma}

\newtheorem{proposition}[thm]{Proposition}
\newtheorem{corollary}[thm]{Corollary}

\newtheorem{conjecture}[thm]{Conjecture}

\renewcommand{\a}{{\mathfrak a}}

\renewcommand{\P}{{\mathbb P}}

\newcommand{\F}{{\mathcal F}}

\newcommand{\J}{{\mathscr J}}
\newcommand{\cK}{{\mathcal K}}
\newcommand{\K}{{\mathbb K}}
\renewcommand{\L}{{\mathscr L}}

\newcommand{\Z}{{\mathbb Z}}

\newcommand{\R}{{\mathbb R}}

\newcommand{\bu}{{\bullet}}
\newcommand{\x}{{\times}}

\newcommand{\la}{{\langle}}
\newcommand{\ra}{{\rangle}}
\newcommand{\lra}{{\,\longrightarrow\,}}

\renewcommand{\th}{{\text{th}}}

\DeclareMathOperator{\card}{cardinality}

\DeclareMathOperator{\im}{im}

\DeclareMathOperator{\Prob}{Prob}
\DeclareMathOperator{\rad}{rad}

\title{Random Variables with Completely Independent Subcollections}
\author{George A. Kirkup}
\address{University of California, Berkeley}
\date{August 16, 2004}
\email{kirkup@math.berkeley.edu}

\begin{document}

\begin{abstract}

We investigate the algebra and geometry of the 
independence conditions on discrete random
variables in which we fix some random
variables and study the complete independence of 
some subcollections.  We interpret such independence conditions
on the random variables as an ideal of algebraic
relations.  After a change of variables, this ideal
is generated by generalized $2 \x 2$ minors of multi-way
tables and linear forms.  In particular, let $\Delta$
be a simplicial complex on some random variables and
$A$ be the table corresponding to the product of those
random variables.  If $A$ is $\Delta$-independent table 
then $A$ can be written as the entrywise sum $A^I + A^0$ 
where $A^I$ is a completely independent table and $A^0$ 
is identically $0$ in its $\Delta$-margins.

We compute the isolated components of the original ideal,
showing that there is only one component that could
correspond to probability distributions, and relate 
the algebra and geometry of the main component to 
that of the Segre embedding.  If $\Delta$ has fewer
than three facets, we are able to compute generators for the 
main component, show that it is Cohen--Macaulay, and give 
a full primary decomposition of the original ideal.
\end{abstract}

\maketitle

\section{Introduction}
\subsection{Set-Theoretic Version of the Main Result}
Let $X_1, \dotsc ,X_n$ be discrete random variables on
the same population.  Then there is an $n$-dimensional
table whose $(i_1,\dotsc,i_n)$ entry is the probability
of $X_j = i_j$ for all $j$.  Given the table of probabilities
$A$ for $X_1,\dotsc,X_n$, and a subset\linebreak $\J \subset \{1,\dotsc,n\}$
it is easy to compute the table, $A_{\J}$ for $X_{j_1},\dotsc,X_{j_r}$ 
by summing over the indices not in $\J$.

The random variables $X_1,\dotsc,X_n$ are called
completely independent if the probabilities satisfy
$$\Prob (X_1 = i_1, \dotsc , X_n = i_n) = \prod_j \Prob(X_j = i_j)$$
for all possible $i_1, \dotsc , i_n$.  If
$\Delta$ is any collection of subsets of $\{1,\dotsc,n\}$
then we say that an $n$-dimensional table is $\Delta$-independent 
if for each $\J \in \Delta$, $A_{\J}$ is completely independent.

With this notation, the main result of this paper implies
\begin{thm}\label{statthm}
If $A$ is a $\Delta$-independent table associated to the 
product variable $X_1 \x \dotsb \x X_n$ then $A$ can be 
written as the (entrywise) sum 
$A^I + A^0$ where $A^I$ is the completely  independent
table and $A^0$ is a table whose margins $A^0_{\J}$ are 
identically $0$ for all $\J \in \Delta$.
\end{thm}

\subsection{The Algebraic Perspective}
The statement above is not the strongest that can 
be made.  Let $A = (x_{i_1,\dotsc,i_n})$ 
be the generic $a_1 \x \dotsb \x a_n$ table over
any field $\K$ and $\Delta$ be any collection of subsets of
$\{1,\dotsc,n\}$.  For each $\J \in \Delta$, $A_{\J}$ is a
table whose entries are sums of the variables
$x_{i_1,\dotsc,i_n}$.  Complete independence of a table, $B$, with
entries in a ring can be expressed by the ideal, $I(B)$, 
generated by generalized $2 \x 2$ minors of the table.
Therefore, $\Delta$-independence of the generic table $A$ 
is expressed by the ideal 
$$I_{\Delta}(A) = \sum_{\J \in \Delta} I(A_{\J}).$$

Theorem \ref{statthm} is implied by a knowledge of  the minimal
primes over $I_{\Delta}$.  We prove that there is only
one minimal prime, $P_{\Delta}$ over $I_{\Delta}$ which does not
contain the sum of all the variables.  Therefore, $P_{\Delta}$
is the only minimal prime that corresponds to probability
distributions.  We parameterize $P_{\Delta}$ and give 
set-theoretic generators for it in terms of the generators 
of a related toric ideal.  In the case in which $\Delta$ has fewer 
than three facets, we compute the generators for 
$P_{\Delta}$ and show that it is a perfect ideal.

The other minimal primes over $I_{\Delta}$ are also
accessible, and we give a fairly complete description
of them.  Moreover, when $\Delta$ has fewer than three facets
we show that $I_{\Delta}$ is a radical ideal.  $I_{\Delta}$
is not always radical and we also give an example in which
$\Delta$ has four facets and $I_{\Delta}$ is not radical.

\subsection{Overview}
In Section \ref{statsec}, we define the principal objects of study
and develop the elementary statistical terminology needed for the
sequel.  Section \ref{covsec} defines the change of variables which
is the foundation for the rest of the exposition.

In Section \ref{minprim} we show that a related toric ideal is contained
in $P_{\Delta}$, and in Section \ref{deltaind} we parameterize $P_{\Delta}$
and give set-theoretic generators for it.  In Section \ref{othersec} we
treat the other minimal primes over $I_{\Delta}$ and show that they
can be understood in terms of $P_{\Delta_i}$ for subcomplexes
$\Delta_i \subset \Delta$.  In Section \ref{prssec} we use
principal radical systems to prove that if $\Delta$ has three or
fewer facets then $P_{\Delta}$ is generated by the set-theoretic
generators given in Section {deltaind} and is a perfect ideal.  We also
prove that in the same case, $I_{\Delta}$ is radical.  Finally, Section
\ref{finalsec} ties up the loose ends with an example in which
$I_{\Delta}$ is not radical, two conjectures and notes on the copmutational
limits encountered.

The main theorems are Theorem \ref{radsegeq} and Theorem \ref{threemain}.
The change of variables in Section \ref{covsec} and the toric
ideal $Q_{\Delta}$ from Section \ref{minprim} are the key technical
points to understand from which Theorem \ref{radsegeq} follows.
Theorem \ref{threemain} is an application of principal radical systems.

\subsection{Acknowledgements}
I was fully supported by the Air Force, through a National 
Defense Science and Engineering Graduate Fellowship.  I would like
to thank my research advisor, David Eisenbud for all his help
and support.  I also received substantial help and encouragement
from Bernd Sturmfels throughout the project.

\section{Statistics for Algebraists}\label{statsec}
\subsection{Random Variables}
A \emph{random variable} $X$ is a function from a set $\Omega$,
a population, to a set $S_X$, the values of $X$.  We define
$$\{X = s\} = X^{-1}(s).$$
If $\Omega$ is finite, we define a new function 
$P_X \colon S_X \lra \R_+$ by
$$P_X(s) = \Prob\{X = s\} = \frac{\card \{X = s\}}{\card \Omega } \cdot$$
$P_X(s)$ can be interpreted as the probability that a randomly 
selected $\omega \in \Omega$ will have $X(\omega) = s$.  A 
\emph{discrete random variable} is a random variable which takes 
finitely many values.  From now on, all our random variables 
will be discrete on a finite population.  That is, $\Omega$ 
and $S_X$ are both finite.

If $X_1, \dotsc , X_n$ are random variables on the same
population, then there is a product variable
$$X_1 \x \dotsb \x X_n \colon \Omega \lra S_{X_1} \x \dotsb \x S_{X_n}$$
defined in the obvious way.  If $X_j$ takes $a_j < \infty$ values,
then there is an $a_1 \x \dotsb \x a_n$ $n$-dimensional (real) table
$$A = (x_{i_1, \dotsc , i_n})$$
whose $(i_1, \dotsc , i_n)$ entry is the probability,
$$\Prob \{X_1 \x \dotsb \x X_n = (i_1, \dotsc , i_n)\}.$$

\subsection{Marginal Tables and Subcollections of Random Variables}\label{plusnot}
Suppose we have an $n$-dimensional array
$A = (x_{i_1, \dotsc , i_n})$ of probabilities associated to 
some random variables $X_1, \dotsc , X_n$.  
Given any $\J = \{j_1, \dotsc ,j_m\} \subset \{1, \dotsc,n\}$
we can define an $a_{j_1} \x \dotsb \x a_{j_m}$ array
which is the probability array for the random variable
$X_{j_1} \x \dotsb \x X_{j_m}$, disregarding the other random
variables.  Such an array is called
an \emph{$m$-margin} of $A$.

To recover the probability of some subcollection of events 
happening, disregarding the other variables, we need 
only to sum over the variables we wish to disregard. For 
example, to disregard the random variable $X_n$, consider
$$\Prob \{X_1 \x \dotsb \x X_{n-1} = (i_1, \dotsc ,i_{n-1})\} = \sum_k x_{i_1, \dotsc ,i_{n-1},k}.$$

In general, suppose that $A = (x_{i_1, \dotsc , i_n})$ is 
an $n$-dimensional array with entries in a ring $R$.
Let $\sigma$ be an ordered $n$-tuple whose $j^\th$ entry,
$\sigma_j$, is either an integer such that $1 \leq \sigma_j \leq a_j$
or the symbol $+$.  Let 
$$\J = \J(\sigma) = \{j_1, \dotsc ,j_m\} = \{j \mid \sigma_j \neq +\}$$
and define
$$x_{\sigma} := \sum_{i_{j} = \sigma_{j} \text{ if } j \in \J} x_{i_1, \dotsc , i_n}.$$
For example, $x_{1,+,3} = \sum_j x_{1,j,3}.$

This essentially allows us to create the desired array, but we need
to index the array correctly.
Fix some $\J = \{j_1, \dotsc ,j_m\} \subset \{1, \dotsc, n\}$
and numbers $i_1, \dotsc , i_m$ such that $1 \leq i_k \leq a_{j_k}$.
We can define a sequence $\sigma(\J)_{\{i_1, \dotsc, i_m\}}$ of length $n$, 
by $\sigma(\J)_k = +$ if $k \not\in \J$, and $\sigma(\J)_{j_k} = i_k$.
Again, let $A = (x_{i_1, \dotsc , i_n})$ be an $n$-dimensional array
with entries in a ring $R$.  We may define an
$a_{j_1} \x \dotsb \x a_{j_m}$ array $A_{\J}$ whose $(i_1, \dotsc ,i_m)$
entry is $x_{(\sigma(\J)_{i_1, \dotsc ,i_m})}.$
This is an $m$-margin of $A$, as described above.

Moreover, if $A$ is an array of probabilities that is associated to 
random variables $X_1, \dotsc , X_n$ and 
$\J \subset \{1, \dotsc ,n\}$, then $A_{\J}$ is the array
of probabilities associated to the random variables
$X_{j_1}, \dotsc , X_{j_m}$.

\subsection{Complete Independence and the Segre Variety}
The random variables $X_1, \dotsc, X_n$ are called 
\emph{completely independent} if the identity
$$\Prob \{X_1 \x \dotsb \x X_n = (i_1, \dotsc ,i_n)\} = \prod_{j=1}^n \Prob \{X_j = i_j\}$$
holds for all values in $S_{X_1} \x \dotsb \x S_{X_n}$.  We
will study the situation in which certain subcollections of the
variables $X_1, \dotsc , X_n$ are completely independent.

Likewise an array $A = (x_{i_1, \dotsc ,i_n})$ with entries in
a ring $R$ will be called \emph{completely independent} if there
are elements of $R$, $\{y_{1,i_1},y_{2,i_2}, \dotsc , y_{n,i_n}\}$,
such that the condition
\begin{equation}\label{statseg}
x_{i_1, \dotsc ,i_n} = \prod y_{j,i_j}
\end{equation}
holds for all choices $(i_1, \dotsc ,i_n)$.

An algebraic geometer will immediately recognize that \eqref{statseg} 
implies that the table $A$ is a point on the Segre variety
\begin{equation}\label{geomsegre}
\P^{a_1 - 1} \x \dotsb \x \P^{a_n-1} \subset \P^{\prod a_j - 1}.
\end{equation}
This brings us to the link between statistics and 
commutative algebra.

\subsection{The Algebraic Definitions}\label{algdef}

The Segre embedding is induced by the ring map
\begin{equation*}
\begin{split}
\sigma : \Z[x_{i_1, \dotsc , i_n}] & \lra \Z[y_{1,i_1},y_{2,i_2}, \dotsc , y_{n,i_n}] \\
   x_{i_1, \dotsc , i_n} & \longmapsto \prod y_{j,i_j}
\end{split}
\end{equation*}
The kernel of $\sigma$, which is the defining ideal of
the Segre variety, can be generated by generalized
$2 \x 2$ minors, which we now define.

As usual, let $A = (x_{i_1, \dotsc , i_n})$ be an $n$-dimensional
array with entries in a ring $R$.
We define a $2 \x 2$ minor about the $l^{\th}$ coordinate of $A$
to be any relation of the form
$$\det \begin{pmatrix} x_{i_1, \dotsc, i_n} & x_{j_1, \dotsc, j_{l-1},i_l,j_{l+1}, \dotsc ,j_n} \\
                       x_{i_1, \dotsc, i_{l-1},j_l,i_{l+1}, \dotsc ,i_n} & x_{j_1, \dotsc, j_n} \end{pmatrix}.$$
This is an interchange of just the $l^{\th}$ coordinate.  Obviously, the ideal 
in $R$ generated by all interchanges of one coordinate will
generate the ideal containing all interchanges of an arbitrary 
number of coordinates.  From \cite[Corollary~1.8]{HTHa}, we 
know that the $2 \x 2$ minors of an $n$-dimensional array 
generate the defining ideal of the Segre embedding.  Thus we
define the \emph{Segre relations} to be these generalized
$2 \x 2$ minors.

We can define an $a_1 \x \dotsb \x a_n$ table with
entries in $R$ to be a map
\begin{equation*}
\begin{split}
B \colon \Z[x_{i_1, \dotsc , i_n}] & \lra R \\
   x_{i_1, \dotsc , i_n} & \longmapsto b_{i_1, \dotsc , i_n}
\end{split}
\end{equation*}
where the $(i_1, \dotsc , i_n)$ entry in $B$ is
defined to be $b_{i_1, \dotsc , i_n}$.  In this language, 
the generic table is the identity map.

We have a diagram
$$\xymatrix{ R \ar[r] & R \otimes_{\Z[x_{i_1, \dotsc , i_n}]} \Z[y_{1,i_1},y_{2,i_2}, \dotsc , y_{n,i_n}]\\
\Z[x_{i_1, \dotsc , i_n}] \ar[r]^-{\sigma} \ar[u]^B & \Z[y_{1,i_1},y_{2,i_2}, \dotsc , y_{n,i_n}] \ar[u]}$$
and we let $I(B) \subset R$ be the kernel of the top map.
This amounts to imposing the Segre relations above on the
table $B$.

Let $A$ be the generic $a_1 \x \dotsb \x a_n$ table and
let $\Delta$ be a collection of subsets of $\{1, \dotsc , n\}$.
Recall the definition of the marginal tables $A_{\J}$ from
Section~\ref{plusnot}.  We define the ideal
$$I_{\Delta} (A) = \sum_{\J \in \Delta} I(A_{\J}).$$
That is, $I_{\Delta}(A)$ is the ideal generated by the
generalized $2 \x 2$ minors of each margin $A_{\J}$, when
$\J \in \Delta$.  We give an example at the end of this section.

This is a special case of what are called ``independence ideals'' 
in the algebraic statistics literature.  See \cite[\S 8.1]{SPE} 
for more about independence models and their corresponding
ideals.  One research paper which gives a discussion of conditional
independence of four random variables is \cite{Matus}.
$I_{\Delta}(A)$ should be thought of as the defining
ideal of the variety of tables which are completely independent
in the margins given by $\Delta$.  We call a table 
\emph{$\Delta$-independent} if it lies on the variety defined 
by $I_{\Delta}(A)$.

If $\J' \subset \J$, then because of the multilinearity of
the Segre relations,  the complete independence of $A_{\J}$ 
implies the complete independence of $A_{\J'}$.  Thus we may
assume that $\Delta$ has the structure of a simplicial complex;
that is, $\J' \subset \J \in \Delta \Longrightarrow \J' \in \Delta$.

The rest of the paper is concerned with the primary decomposition
of the ideals $I_{\Delta}(A)$.  For any $\Delta$ we will show there is only
one minimal prime which does not contain $x_{+, \dotsc , +}$.  This
component is the most important because when $A$ represents a 
probability distribution, $x_{+, \dotsc , +} = 1$.  Thus we study 
that prime and relate it algebraically and geometrically to the 
Segre variety.  When $\Delta$ is a simplicial complex with three
or fewer facets, we can compute generators for the main component
and show that it is perfect.  In that case we will also show that
that $I_{\Delta}(A)$ is a radical ideal and give a full
primary decomposition.

Throughout the exposition, we will consider the following running
example for clarity:  $n=3$, $a_1=a_2=a_3=2$, and
$$\Delta = \{\{1,2\},\{1,3\},\{2,3\},\{1\},\{2\},\{3\},\emptyset\}.$$
In this case $R = \K[x_{i,j,k}]$ is a polynomial ring with $8$
variables and $I_{\Delta}$ is generated by $3$ elements:
$$I_{\Delta} = \la \det \begin{pmatrix} x_{1,1,+} & x_{1,2,+} \\ x_{2,1,+} & x_{2,2,+} \end{pmatrix},\det \begin{pmatrix} x_{1,+,1} & x_{1,+,2} \\ x_{2,+,1} & x_{2,+,2} \end{pmatrix},\det \begin{pmatrix} x_{+,1,1} & x_{+,1,2} \\ x_{+,2,1} & x_{+,2,2} \end{pmatrix} \ra.$$
Despite its appearance, $I_{\Delta}$ is not a binomial ideal because
$x_{1,1,+} = x_{1,1,1}+x_{1,1,2}$.

\section{A Linear Change of Variables}\label{covsec}
\subsection{Set-theoretic Heuristics}
Let $\Delta$ be some fixed collection of subsets of $\{1,\dotsc,n\}$.  
Our goal is to decompose the ideal 
$I_{\Delta}(A) \subset R = \K[x_{i_1, \dotsc , i_n}]$
which is defined by the complete independence of the collection
of margins of the generic table $A$ given by $\Delta$.  First,
it will be helpful and illuminating to perform a linear change
of variables on $R$ which makes $I_{\Delta}$ an ideal generated
by quadratic binomials and linear forms.  We will show 
that $R/I_{\Delta}(A)$ is a polynomial ring over a ring of 
smaller dimension.

Set-theoretically, suppose that one table $A$ is $\Delta$-independent, 
and another table $B$ has the property that for each $\J \in \Delta$, 
$B_{\J} = 0$.
Then the sum (entry by entry) $A + B$ is also $\Delta$-independent.
This is a trivial result of the fact that the equations which define
$\Delta$-independence only involve entries of the marginal tables and
$B$ is identically $0$ in its $\Delta$-margins.  In this section, 
we will develop this idea algebraically.

\subsection{$S_{\Delta}$, $T_{\Delta}$ and the Change of Variables}\label{TStau}
We define a ring $S_\Delta$, to be the polynomial ring over $\K$ with variables
that are indexed by the entries in the marginal tables given by the
elements of $\Delta$.
That is, for every $\J \in \Delta$, $A_{\J} = (x_{i_1, \dotsc , i_n})$
with $i_k = +$ for every $k \not\in \J$.  So for every $\J \in \Delta$, create 
a formal symbol $X_{i_1, \dotsc , i_n}$ with $i_k = \bu$ for every $k \not\in \J$ 
and  $1 \leq i_j \leq a_j$ for all $j \in \J$.  Then let $S_\Delta$ be the
polynomial ring over $\K$ generated by these formal symbols.

Now consider the map of rings $\tau_\Delta \colon S_\Delta \lra R$ defined by 
$$X_{i_1, \dotsc , i_n} \longmapsto x_{i_1, \dotsc , i_n}$$
in which $\bu$ changes to $+$.
The kernel of $\tau_\Delta$, $K_\Delta \subset S_\Delta$, is generated by linear
forms.  Let $T_\Delta = S_\Delta/K_\Delta$ be the coordinate ring of 
\emph{$\Delta$-marginal tables}.  Set-theoretically, a 
$\Delta$-marginal table $B$ represents the class of tables $B'$ such
that for all $\J\in \Delta$, $B_\J = {B'}_{\J}$.

If $\Delta$ and $\Delta'$ have the property
that the maximal elements of $\Delta$ and $\Delta'$ are the same,
it is clear that $T_{\Delta} \cong T_{\Delta'}$.  Since there is no 
ambiguity in $T_{\Delta}$, we will replace all $\bu$'s in the indices 
of the variables by $+$'s as usual.

On the other hand, let 
$$Z_\Delta = R/(L_{\Delta}(A)),$$
the coordinate ring of tables whose margins given by $\Delta$ are
identically zero.  Since the ideal we quotient by is generated
by linear forms, $Z_\Delta$ is a polynomial ring over $\K$.
Moreover, since the image of $\tau_\Delta$ is generated by the
linear forms which generate $\sum_{\L \in \Delta} L_{\L}$,
we have 
\begin{equation}\label{vecbund}
T_\Delta \otimes Z_\Delta \cong \im \tau_\Delta \otimes Z_\Delta \cong R.
\end{equation}
Set-theoretically, this says that the space of $a_1 \x \dotsb \x a_n$
tables is a trivial bundle over the space of $\Delta$-marginal tables.

\begin{proposition}\label{cov}
Suppose that $I = \la f_n \ra$ is any ideal in $R$
such that the $f_n$ are written entirely in terms of the margins
given by $\Delta$, as above.  Then let $F_n$ be the polynomial in 
$S_\Delta$ (or $T_\Delta$) which has the same form as $f_n$ except that 
the lower-case $x$'s are replaced by upper-case $X$'s and the
$+$'s are replaced by $\bu$'s.  Let
$I^{S_\Delta} := K_\Delta + \la F_n \ra$.

Then $I$ is prime (respectively radical, perfect) if and
only if $I^{S_\Delta}$ is prime (respectively radical, perfect).
Moreover, the Betti diagram of $I$ as an $R$-module is the
same as that of $I^{T_{\Delta}}$ as a $T_{\Delta}$-module.
\end{proposition}
\begin{proof}
Since polynomial rings are flat over the ground field, by \eqref{vecbund}
$$R/I \cong S_\Delta/(I^{S_\Delta}) \otimes Z_\Delta,$$
which is a polynomial ring over $S_\Delta/(I^{S_\Delta})$.
Thus, $R/I$ is a domain (resp. reduced, Cohen--Macaulay) if and
only if $S_\Delta/I^{S_\Delta}$ is a domain (resp. reduced, Cohen--Macaulay).
\end{proof}

\subsection{Generators for $K_{\Delta}$ and Our Example}
We can also describe the generators of $K_\Delta$.  The idea is that if
we have two margins $A_{\J}$ and $A_{\cK}$ then they have
an ``intersection'' which is $A_{\J \cap \cK}$.  In 
particular, the entries of $A_{\J \cap \cK}$ will have
a representation as sums of elements of $A_{\J}$ and
$A_{\cK}$, and they must agree.  For ease of notation, we will
assume that $\J = \{1, \dotsc, r\}$ and $\cK = \{s, \dotsc ,n\}$,
so $\L = \{r, \dotsc ,s\}$.  Then we have an ideal of relations
$$R_{\J, \cK} = \la \sum_{i_m \mid m < s} X_{i_1, \dotsc , i_r, +, \dotsc, +} - 
\sum_{i_m \mid m > r} X_{+, \dotsc, +, i_s, \dotsc, i_n} \ra$$
for all choices of $(i_s, \dotsc , i_r)$.  

\begin{proposition}
$K_\Delta$ is generated by $\sum R_{\J, \cK}$ for all pairs of 
$\J, \cK \in \Delta$.
\end{proposition}

We now turn to our example, in which $n=3$, and
$$\Delta = \{\{1,2\},\{1,3\},\{2,3\},\{1\},\{2\},\{3\},\emptyset\}.$$
In this case,
$$S_\Delta = \K[X_{i,j,\bu},X_{i,\bu,k},X_{\bu,j,k},X_{i,\bu,\bu},X_{\bu,j,\bu},X_{\bu,\bu,k},X_{\bu,\bu,\bu}]$$
for $1 \leq i,j,k \leq 2$.  Thus $S_{\Delta}$ has 19 variables.
$K_{\Delta}$ is easy to describe.  For example,
$$R_{\{1,2\}, \{1,3\}} = \la (X_{1,1,\bu}+X_{1,2,\bu}) - (X_{1,\bu,1}+X_{1,\bu,2}), (X_{2,1,\bu}+X_{2,2,\bu}) - (X_{2,\bu,1}+X_{2,\bu,2}) \ra,$$
$$R_{\{1,2\},\{1\}} = \la (X_{1,1,\bu}+X_{1,2,\bu}) - X_{1,\bu,\bu}, (X_{2,1,\bu}+X_{2,2,\bu}) - X_{2,\bu,\bu} \ra,$$
$$R_{\{1,2\},\{3\}} = \la (X_{1,1,\bu}+X_{1,2,\bu} + X_{2,1,\bu}+X_{2,2,\bu}) - (X_{\bu,\bu,1}+X_{\bu,\bu,2}) \ra.$$
$K_{\Delta}$ can be generated by $30$ linear forms, but of
course this is not minimal, as illustrated by the inclusion
$$R_{\{1,3\},\{1\}} \subset R_{\{1,2\}, \{1,3\}} + R_{\{1,2\},\{1\}}.$$
$K_{\Delta}$ can be minimally generated by $12$ linear forms
so $T_{\Delta} = S_{\Delta}/K_{\Delta}$ is a polynomial ring
of dimension $7$.

Notice that if $\Delta' = \{\{1,2\},\{1,3\},\{2,3\}\}$, $S_{\Delta'}$
has $12$ variables, and $K_{\Delta'}$ is minimally generated by $5$
linear forms, and $T_{\Delta'} \cong T_{\Delta}$.

\section{In Search of a Statistically Significant Component $I_{\Delta} (A)$}\label{minprim}
\subsection{A Related Toric Ideal}
In the following sections we will prove that there is only one
minimal prime over $I_{\Delta} (A)$, for a generic $a_1 \x \dotsb \x a_n$ 
table $A$, which does not contain $x_{+,\dotsc,+}$.  This will
be the only statistically significant component of $I_{\Delta}$
because when $A$ is a probability distribution, $x_{+,\dotsc,+}=1$.
We will identify the main component as the kernel of ring map, and 
relate it to a toric ideal.

The first step is to define the toric ideal.  Let $\Delta$ be a collection 
of subsets of $\{1,\dotsc,n\}$ and let
\begin{equation*}
\begin{split}
\eta_\Delta \colon S_\Delta \lra &\K[y_{i,j_i} \mid 1\leq j_i \leq a_i \text{ or } j_i = \bu] \\
  X_{j_1,\dotsc,j_n} \longmapsto &\prod y_{i,j_i}
\end{split}
\end{equation*}
Finally, let $Q_{\Delta} = \ker \eta_{\Delta}$.  Since $Q_{\Delta}$
is defined as the kernel of a monomial map, it is generated by binomials.
The rest of this section will be devoted to showing that $Q_{\Delta}$ 
is contained in $(I_{\Delta} : x_{+,\dotsc,+}^{\infty})$.  

In our example, where $\Delta$ has facets $\{1,2\},\{1,3\},\{2,3\}$,
$Q_{\Delta}$ is generated by such binomials as 
$$\det \begin{pmatrix} X_{1,1,\bu} &  X_{2,1,\bu} \\ X_{1,\bu,1} & X_{2,\bu,1} \end{pmatrix} \text{ and } \det \begin{pmatrix}X_{1,1,\bu}& X_{2,1,\bu} \\ X_{1,\bu,\bu} & X_{2,\bu,\bu} \end{pmatrix}.$$
$S_{\Delta}/Q_{\Delta}$ is a ring of dimension $6$.

\subsection{Some Useful Elements of the Ideal $I(A)$}
First we will construct elements in $I(A)$ which
will allow us to view $+$ like any other index.

\begin{proposition}\label{summin}
Let $A$ be the generic $a_1 \x \dotsc \x a_n$ table.  Then
$$\det \begin{pmatrix} x_{i_1, \dotsc, i_n} & x_{j_1, \dotsc, j_{l-1},i_l,j_{l+1}, \dotsc ,j_n} \\
x_{i_1, \dotsc, i_{l-1},+,i_{l+1}, \dotsc ,i_n} & x_{j_1, \dotsc,j_{l-1},+,j_{l+1}, \dotsc ,j_n} \end{pmatrix} \in I(A).$$
\end{proposition}
\begin{proof}
Consider the sum
$$\sum_{k} \det \begin{pmatrix} x_{i_1, \dotsc, i_{l-1},i_l,i_{l+1}, \dotsc ,i_n} & x_{j_1, \dotsc, j_{l-1},i_l,j_{l+1}, \dotsc ,j_n} \\
   x_{i_1, \dotsc, i_{l-1},k,i_{l+1}, \dotsc ,i_n} & x_{j_1, \dotsc,j_{l-1},k,j_{l+1}, \dotsc ,j_n} \end{pmatrix}$$
which by the multilinearity of the minors is
$$\det \begin{pmatrix} x_{i_1, \dotsc, i_{l-1},i_l,i_{l+1}, \dotsc ,i_n} & x_{j_1, \dotsc, j_{l-1},i_l,j_{l+1}, \dotsc ,j_n} \\
 \sum_k x_{i_1, \dotsc, i_{l-1},k,i_{l+1}, \dotsc ,i_n} & \sum_k x_{j_1, \dotsc,j_{l-1},k,j_{l+1}, \dotsc ,j_n} \end{pmatrix}.$$
By the definition of $+$ notation from Section~\ref{plusnot}, this is
$$\det \begin{pmatrix} x_{i_1, \dotsc, i_n} & x_{j_1, \dotsc, j_{l-1},i_l,j_{l+1}, \dotsc ,j_n} \\
x_{i_1, \dotsc, i_{l-1},+,i_{l+1}, \dotsc ,i_n} & x_{j_1, \dotsc,j_{l-1},+,j_{l+1}, \dotsc ,j_n} \end{pmatrix}$$
which establishes the result.
\end{proof}

This proposition allows us to let any number of coordinates 
equal ``+'', and interchange them freely.

As an example, consider the case in which $n = 2$.  For any $i_1,i_2$
$$x_{i_1,i_2} x_{+,+} - x_{i_1,+}x_{+,i_2} \in I_{\{1,2\}}(A).$$
If the $x_{i,j}$ are really probabilities, then $x_{+,+} = 1$ so this
relation becomes $x_{i_1,i_2} = x_{i_1,+}x_{+,i_2}$, which is the
independence condition for two random variables, as in \eqref{statseg}.

\subsection{An Intermediate Ideal, $J_{\Delta} \subset Q_{\Delta}$}
There are some quadratic binomials in $Q_{\Delta}$
which play a special role in the discussion.  Let 
$J_{\Delta} \subset Q_{\Delta}$ be generated by binomials
$$f = X_{\bar{\imath}_1} \cdot X_{\bar{\imath}_2} - X_{\bar{\jmath}_1} \cdot X_{\bar{\jmath}_2} \in Q_{\Delta}$$
such that $X_{\bar{\imath}_1},X_{\bar{\jmath}_1}$ are both entries in $A_{\J}$
for some $\J \in \Delta$.  Since $f \in Q_{\Delta}$, this implies that
$X_{\bar{\imath}_2},X_{\bar{\jmath}_2}$ are both entries in $A_{\cK}$ for some
$\cK \in \Delta$.

In our example, $J_{\Delta}$ will be generated by $I_{\Delta}$ and the
$2\x 2$ minors of the three matrices symmetric to
\begin{equation}\label{flatmat}
\begin{pmatrix} X_{1,1,\bu} & X_{1,2,\bu} & X_{1,\bu,\bu} & X_{1,\bu,1} & X_{1,\bu,2} \\
                X_{2,1,\bu} & X_{2,2,\bu} & X_{2,\bu,\bu} & X_{2,\bu,1} & X_{2,\bu,2} \end{pmatrix}.
\end{equation}

\begin{lemma}
Let $f = X_{\bar{\imath}_1} \cdot X_{\bar{\imath}_2} - X_{\bar{\jmath}_1} \cdot X_{\bar{\jmath}_2}$
be a generator of $J_{\Delta}$ such that $X_{\bar{\imath}_1}$ is an
entry in $A_{\J}$ and $X_{\bar{\imath}_2}$ is an entry in $A_{\cK}$.
Then
$$L_{\J \cap \cK} \cdot \la f \ra \subset I_{\Delta}.$$
\end{lemma}
\begin{proof}
The proof is very technical (but elementary).  In our running
example, the result follows from the following line of reasoning.  
The matrix \eqref{flatmat} has the property that the first
$3$ columns and the last $3$ columns have rank $1$.  Since they
share the middle column, either each column of \eqref{flatmat}
is a scalar multiple of the middle column, or the middle
column is identically $0$.  Thus, either the $2 \x 2$ minors
of \eqref{flatmat} vanish or $X_{1,\bu,\bu} = X_{2,\bu,\bu}=0$.

Now we turn to the detailed proof.
Since all the calculations will happen in the margin 
$A_{\J \cup \cK}$, we can assume that 
$\{1, \dotsc , n\} = \J \cup \cK$ for ease of notation.  
We re-index so that $\J = \{1, \dotsc , s\}$ and
$\cK = \{r, \dotsc , n\}$, so $\L = \{r, \dotsc , s\}$.
After this reorganization, $f$ is the following determinant
$$q = \det \begin{pmatrix} X_{i_{1,1},\dotsc ,i_{1,s},+,\dotsc, +} & X_{+,\dotsc ,+,j_{2,r}, \dotsc, j_{2,n}} \\
             X_{j_{1,1},\dotsc ,j_{1,s},+,\dotsc, +} & X_{+,\dotsc ,+,i_{2,r}, \dotsc, i_{2,n}} \end{pmatrix}$$
where $i_{1,k} = j_{1,k}$ for all $k < r$ and 
$i_{2,k} = j_{2,k}$ for all $k >s$.  Moreover,
$\{i_{1,k},i_{2,k}\} = \{j_{1,k},j_{2,k}\}$ for each
$r \leq k \leq s$.  Thus, $f$ can be thought of as the
exchange of some number of indices between
$X_{\bar{\imath}_1}$ and $X_{\bar{\imath}_2}$.  Clearly, these
exchanges can be generated by exchanges of one coordinate.
Re-index again, so that $f$ is an exchange of the $r^{\th}$
coordinate.  Then $f$ can be written as
$$f = \det \begin{pmatrix} X_{j_1,\dotsc,j_{r-1},j_r,j_{r+1},\dotsc,j_s,+,\dotsc,+} & X_{+,\dotsc ,+,k_r,j_{r+1},\dotsc ,j_n} \\
X_{j_1,\dotsc,j_{r-1},k_r,j_{r+1},\dotsc,j_s,+,\dotsc,+} & X_{+, \dotsc, +,j_r,k_{r+1} \dotsc ,k_n} \end{pmatrix}.$$

If $l = x_{+, \dotsc ,+,i_r, \dotsc ,i_s,+, \dotsc ,+}$
is any generator of $L_{\L}$, then we need to show that $lf$ 
is in $I(A_{\J}) + I(A_{\cK})$.  We will constuct this
product explicitly.

Consider the sum
\begin{equation*}
\begin{split}
x_{j_1, \dotsc ,j_{r-1}, i_r,j_{r+1} \dotsc ,j_s,+,\dotsc ,+} 
&\det \begin{pmatrix} 
x_{+,\dotsc,+,k_r,k_{r+1},\dotsc,k_n}&x_{+,\dotsc,+,k_r,i_{r+1},\dotsc,i_s,+,\dotsc,+} \\
x_{+,\dotsc,+,j_r,k_{r+1},\dotsc,k_n}&x_{+,\dotsc,+,j_r,i_{r+1},\dotsc,i_s,+,\dotsc,+}
\end{pmatrix} +\\
x_{+,\dotsc,+,k_r,k_{r+1},\dotsc,k_n}
&\det \begin{pmatrix} 
x_{j_1,\dotsc,j_{r-1},j_r,j_{r+1},\dotsc,j_s,+,\dotsc,+} & x_{+,\dotsc,+,j_r,i_{r+1},\dotsc,i_s,+,\dotsc,+} \\
x_{j_1,\dotsc,j_{r-1},i_r,j_{r+1},\dotsc,j_s,+,\dotsc,+} & x_{+,\dotsc,+,i_r,i_{r+1},\dotsc ,i_s,+,\dotsc,+}
\end{pmatrix} + \\
x_{+,\dotsc,+,j_r,\dotsc,j_s,k_{s+1},\dotsc ,k_n}
&\det \begin{pmatrix} 
x_{j_1,\dotsc,j_{r-1},i_r,j_{r+1},\dotsc,j_s,+,\dotsc ,+} & x_{+,\dotsc,+,i_r,i_{r+1},\dotsc,i_s,+,\dotsc,+}\\
x_{j_1,\dotsc,j_{r-1},k_r,j_{r+1},\dotsc,j_s,+,\dotsc ,+} & x_{+,\dotsc,+,k_r,i_{r+1},\dotsc,i_s,+,\dotsc,+}
\end{pmatrix}
\end{split}
\end{equation*}
which is also evidently equal to
$$(X_{+, \dotsc ,+,i_r, \dotsc ,i_s,+, \dotsc ,+}) f = lf \in I(A_{\J}) + I(A_{\cK}).$$
This completes the calculation.
\end{proof}

In our example, in which $\Delta$ has facets
$\{\{1,2\},\{1,3\},\{2,3\}\}$, $J_{\Delta} = Q_{\Delta}$.
This is a result of the fact that $\Delta$ has three facets.
The smallest example in which $J_{\Delta}\neq Q_{\Delta}$
is when $\Delta$ has facets $$\{\{1,2\},\{1,3\},\{2,4\},\{3,4\}\}.$$
In this case,
$$X_{1,1,+,+}X_{+,+,1,1}-X_{1,+,1,+}X_{+,1,+,1}$$
is in $Q_{\Delta}$ but not $J_{\Delta}$.

\subsection{The Relationship Between $Q_{\Delta}$ and $I_{\Delta}$}

We are ready for the main result of the section.

\begin{proposition}\label{Qcolon}
Let $\Delta$ be any simplicial complex, and let
$\L$ be the intersection of the facets of $\Delta$.
If $f$ is any binomial of degree $s$ in $Q_{\Delta}$
then
$$f \in (I_{\Delta} : L_{\L}^s).$$
In particular, $Q_{\Delta} \subset (I_{\Delta}:X_{+,\dotsc,+}^{\infty})$.
\end{proposition}
\begin{proof}

The proof of this proposition is, again, rather technical.
The main obstacle is notation which easily gets confusing.
Thus, we here give an example which will serve to show how
to prove the general result.

As our example, we choose the case in which $\Delta$ has
facets 
$$\{\{1,2\},\{1,3\},\{2,4\},\{3,4\}\},$$ and let
$$f = X_{1,1,+,+}X_{+,+,1,1}-X_{1,+,1,+}X_{+,1,+,1}.$$
We will show that $X_{+,+,+,+}^2 f \in I_{\Delta}$.
Since $J_{\Delta} \subset (I_{\Delta} : X_{+,+,+,+}),$
it suffices to show that $X_{+,+,+,+} f \in J_{\Delta}$.
\begin{equation*}
\begin{split}
&X_{+,+,+,+}(X_{1,1,+,+}X_{+,+,1,1}-X_{1,+,1,+}X_{+,1,+,1}) = \\
&= X_{+,+,+,+} X_{1,1,+,+} X_{+,+,1,1} - X_{+,+,+,+}X_{1,+,1,+}X_{+,1,+,1}\\
&= X_{1,+,+,+} X_{+,1,+,+} X_{+,+,1,1} - X_{1,+,+,+}X_{+,+,1,+}X_{+,1,+,1}\\
&= X_{1,+,+,+} (X_{+,1,+,+} X_{+,+,1,1} - X_{+,+,1,+}X_{+,1,+,1})
\end{split}
\end{equation*}
Since $X_{+,1,+,+} X_{+,+,1,1} - X_{+,+,1,+}X_{+,1,+,1}$ is
in $J_{\Delta}$, the example is shown.
\end{proof}

\section{$\Delta$-Independence and Complete Independence}\label{deltaind}

\subsection{The Segre Embedding, $\sigma_{\Delta}$, and $P_{\Delta}$}
In this section we study the relationship between tables
which are $\Delta$-independent and tables which are completely 
independent.  It is obvious that any table which is competely
independent is also $\Delta$-independent.  By Proposition~\ref{cov},
we know that inside the variety of $\Delta$-independent tables
is a trivial bundle over the Segre variety.  We will extablish a 
close connection between the ideal $K_{\Delta}+Q_{\Delta}$
and the defining ideal of the Segre variety.  In this section, we
will assume that $\Delta$ is a simplicial complex.

The variety of completely independent tables,
or the Segre variety, can be parameterized by
\begin{equation*}
\begin{split}
\sigma_{\{1,\dotsc,n\}} \colon R = \K[x_{i_1, \dotsc , i_n}] & \lra \K[y_{1,i_1}, \dotsc , y_{n,i_n}] \\
                                   x_{i_1, \dotsc , i_n}     & \longmapsto \prod_j y_{j,i_j}
\end{split}
\end{equation*}
as in Section~\ref{algdef}.
This map may be composed with $\tau_\Delta$ from Section 
\ref{TStau} to give a map
$$\sigma_\Delta \colon S_\Delta \lra \K[y_{1,i_1}, \dotsc , y_{n,i_n}].$$
Let $P_{\Delta}$ be the kernel of $\sigma_\Delta$.
Thus $P_{\Delta}$ is a prime ideal which defines
the variety of $\Delta$-marginal tables which come from
a point on the Segre variety.

We Have the following commutative diagram,
$$\xymatrix{ & R = \K[x_{i_1,\dotsc,i_n}] \ar[rd]^{\sigma_{\{1,\dotsc,n\}}} & \\
S_{\Delta} \ar[ru]^{\tau_{\Delta}}\ar[rr]^{\sigma_{\Delta}}\ar[rd]^{\eta_{\Delta}}&&\K[y_{1,i_1}, \dotsc , y_{n,i_n}] \\
 & \K[y_{1,i_1}, \dotsc , y_{n,i_n} \mid 1 \leq i_j \leq a_i \text{ or } i_j = \bu] \ar[ru]^{\sum_{\bu}} &}$$

\subsection{The Main Theorem}

We are ready for the main theorem, of
which Theorem \ref{statthm} is a corollary.
The above commutative diagram summarizes all
the main definitions.

\begin{thm}\label{radsegeq}
If $\Delta$ is any simplicial complex,
$$P_{\Delta} = \rad (K_{\Delta} + Q_{\Delta})$$
where $P_{\Delta} = \ker \sigma_{\Delta}$,
$K_{\Delta} = \ker \tau_{\Delta}$ and
$Q_{\Delta} = \ker \eta_{\Delta}$.
\end{thm}
\begin{proof}
First, we need to show that 
$$Q_\Delta + K_{\Delta} \subset P_\Delta.$$
It suffices to show that 
$\sigma_{\Delta} (Q_{\Delta}+K_{\Delta}) = 0$,  
which is clear by the definitions.

On the other hand, let $B$ be any point in $S_{\Delta}$ on
$V(K_{\Delta}+Q_{\Delta})$.  Since it is a 
point on $V(Q_{\Delta})$, it can be represented by
$b_{i_1,\dotsc,i_n} = \prod_j y_{j,i_j}$.
Suppose that $\J \in \Delta$ such that
$B_{\J} \neq 0$.  Then re-index so that
$\J = \{1,\dotsc,r\}$ and $b_{1,\dotsc,1,\bu,\dotsc,\bu} \neq 0$.
Now take any $j \in \J$ (and re-index so $j=1$).
Since $B$ is a point on $V(K_{\Delta})$,
\begin{equation*}
\begin{split}
y_{1,\bu} \prod_{2}^r y_{j,1} \prod_{r+1}^n y_{j,\bu} &= b_{\bu,1,\dotsc,1,\bu,\dotsc,\bu} \\
          &= \sum_1^{a_1} b_{i,1, \dotsc,1,\bu,\dotsc,\bu} \\
          &= y_{1,+} \prod_{2}^r y_{j,1} \prod_{r+1}^n y_{j,\bu}.
\end{split}
\end{equation*}
Since $\prod_{2}^r y_{j,1} \prod_{r+1}^n y_{j,\bu} \neq 0$, that
means that $y_{1,\bu} = y_{1,+}$.  Therefore, if $j$ is any 
index such that there is a face $j \in \J \in \Delta$ with 
$B_{\J} \neq 0$ then $y_{j,\bu} = y_{j,+}$.

Therefore, if for each $j$ there is a $\J \in \Delta$
such that $B_{\J} \neq 0$. Then we can let 
$z_{j,i_j} = y_{j,i_j}$ and
$$b_{i_1,\dotsc,i_n} = \prod y_{j,i_j}$$
since $y_{j,\bu} = y_{j,+}$ for all $j$.  Therefore,
$B$ in in $V(P_{\Delta})$.

On the other hand suppose that $\cK$ is the
maximal set such that for any face $\J \in \Delta$,
if $\J \cap \cK \neq \emptyset$, $B_{\J}=0$.
Re-index so that $\cK = \{1,\dotsc,r\}$.  If there
is any face $\J$ such that $B_{\J} \neq 0$,
$\J$ must be disjoint from $\cK$.  Re-index again
so that face is $\{r+1,\dotsc,s\}$, and we have
$$\prod_1^r y_{j,\bu} \cdot \prod_{r+1}^s y_{j,i_j} \cdot \prod_{s+1}^n y_{j,\bu} \neq 0.$$
Therefore, $y_{j,\bu} \neq 0$ for any $j\leq r$.

If $\L \in \Delta$ such that $\L \cap \cK = \emptyset$
then let $\L' = \L \smallsetminus \cK$.  Since
$B_{\L} = 0$, $B_{\L'} = 0$ also.  Re-index so that
$\L' = \{r+1,\dotsc,s\}$.  Then
$$\prod_1^r y_{j,\bu} \cdot \prod_{r+1}^s y_{j,i_j} \cdot \prod_{s+1}^n y_{j,\bu} = 0$$
for all $y_{j,i_j}$, $r < j \leq s$.  Therefore,
either there is some $r < j \leq s$ such that 
$y_{j,i_j} = 0$ or there is some $j > s$ such
that $y_{j,\bu} = 0$.  The former case is impossible
since that would mean $B_{\J} = 0$ for any
$\J \ni j$ which contradicts the maximality of
$\cK$.  Therefore, for any $\L$ which intersects
$\cK$, there is some $j \not\in \cK \cup \L$ such
that $y_{j,\bu} = 0$.  Since $j \not\in \cK$ we know
that $y_{j,+} = y_{j,\bu}$.

Therefore, let 
\begin{eqnarray*}
  z_{j,1} = &y_{j,\bu} & \;\;\; \text{ for $j \in \cK$} \\
z_{j,i_j} = &0         & \;\;\; \text{ for $j \in \cK$, $i_j > 1$} \\
z_{j,i_j} = &y_{j,i_j} & \;\;\; \text{ for $j \not\in \cK$}. \\
\end{eqnarray*}
Notice that $z_{j,+} = y_{j,\bu}$ for all $j$.  By the
previous paragraph, if $\J \cap \cK \neq \emptyset$
then $B_{\J} = 0$.  Moreover, if $\J \cap \cK = \emptyset$
and $b_{i_1,\dotsc,i_n}$ is a coordinate of $\J$ then
$$\prod z_{j,i_j} = \prod_{j\in \J} y_{j,i_j} \prod_{j \not\in \J} y_{j,\bu} = b_{i_1,\dotsc,i_n}$$
so $B \in V(P_{\Delta})$

We have thus shown that 
$V(K_{\Delta}+Q_{\Delta}) = V(P_{\Delta})$,
which implies that 
$\rad(K_{\Delta}+Q_{\Delta}) = P_{\Delta}$
since $P_{\Delta}$ is prime.
\end{proof}

\begin{corollary}
$P_{\Delta}$ is the only minimal prime over $I_{\Delta}$ which
does not contain $x_{+,\dotsc,+}$.
\end{corollary}
\begin{proof}
By Theorem \ref{radsegeq}, $P_{\Delta}$ is the only minimal
prime over $K_{\Delta}+Q_{\Delta}$.  Therefore, by
Proposition~\ref{Qcolon}, $P_{\Delta}$ is
the only minimal prime over $I_{\Delta}$ 
which does not contain $x_{+,\dotsc,+}$.
\end{proof}

Now we state Theorem \ref{radsegeq} in a set-theoretic
form, which slightly generalizes Theorem \ref{statthm}.
\begin{corollary}\label{statcor}
Let $\K$ be any field and $B$ be any table with
entries in $\K$ which is $\Delta$-independent.
If the sum of the entries in $B$ is not $0$,
then $B$ can be written as the (entrywise) sum 
$B^I + B^0$ where $B^I$ is the completely independent 
table whose $(i_1,\dotsc,i_n)$ entry is
$$\prod_j B_{+,\dotsc,+,i_j,+,\dotsc,+}$$
and $B^0$ is a table whose $\Delta$-margins are 
identically $0$.
\end{corollary}

\subsection{An Application to Computational Statistics}
Suppose $B$ is any $a_1 \x \cdots \x a_n$
table, which is the probability distribution for
a random variable $X_1 \x \cdots \x X_n$
and we want to know which sets of the random
variables are completely independent.

If $B$ is any table, let $B^I$ be the table
whose $i_1,\dotsc,i_n$ entry is
$$(B^I)_{i_1,\dotsc,i_n} = \prod_j B_{+,\dotsc,+,i_j,+,\dotsc,+}$$
and let $B^0 = B-B^I$, the entrywise difference
of $B$ and $B^I$.

There is a simplicial complex $\Delta(B)$ such
that $\J \in \Delta(B)$ if and only if
$(B^0)_{\J} = 0$.  Therefore, by Corollary
\ref{statcor}, $\Delta(B)$ gives exactly 
the collection of subsets of
$\{X_1,\dotsc,X_n\}$ which are completely
independent.

This gives an algorithm for determining which
collections of the random variables are
completely independent which is more efficient
than the obvious one which computes each margin 
and determines if it lies on the Segre variety.

\section{The Other Minimal Primes Over $I_{\Delta}(A)$}\label{othersec}

\subsection{Some Technical Results}
Having extablished that $P_{\Delta}$ is the only minimal prime over
$I_{\Delta}(A)$ not containing $x_{+,\dotsc,+}$, it remains to discuss
the minimal primes over $I_{\Delta}(A)$ which do contain $x_{+,\dotsc,+}$.
The follwing simple, technical result, which explains the interplay 
between the $L_{\L}$ and $I(A_{\cK})$, will be the foundation of the
discussion.

\begin{proposition}\label{LIint}
Suppose that $\L,\cK, \J_1 , \J_2$
are subsets of $\{1, \dotsc ,n\}$ such that
$$\cK \supset \J_1 \cup \J_2,$$
$$\L \supset \J_1 \cap \J_2 \cap \cK.$$
Then
$$L_{\J_1} \cdot L_{\J_2} \subset L_{\L} + I(A_{\cK}).$$
\end{proposition}
\begin{proof}
It is clear that if $\L \subset \L '$ then
$L_{\L} \subset L_{\L '}$
so we may assume that 
$$\cK \supset \L = \J_1 \cap \J_2 \cap \cK.$$

Since all the calculations will be done in $A_{\cK}$, we will
assume that $\cK = \{1, \dotsc, n\}$.  Then we re-index so that
$\J_1 =  \{1, \dotsc ,s\}$ and $\J_2 = \{r, \dotsc , n\}$ so
$\L = \{r, \dotsc , s\}$.

Let $x_{b_1, \dotsc, b_s, +, \dotsc , +}$ and
$x_{+,\dotsc,+,c_r, \dotsc, c_n}$ be arbitrary generators of
$L_{\J_1}$ and  $L_{\J_2}$ respectively.
Now consider the following element of $I(A)$.
$$\det\begin{pmatrix} x_{b_1, \dotsc, b_{r-1},c_r, \dotsc, c_n} & x_{b_1, \dotsc ,b_s, +, \dotsc ,+} \\
       x_{+,\dotsc,+,c_r, \dotsc, c_n} & x_{+,\dotsc,+,b_r, \dotsc, b_s, +, \dotsc ,+} \end{pmatrix}$$
The result is clear since $x_{+,\dotsc,+,b_r,\dotsc,b_s,+,\dotsc,+} \in L_{\L}$.
\end{proof}

\begin{corollary}\label{LIintmp}
Suppose that $\L \subsetneq \cK$ and $Q$ is
a prime ideal containing $L_{\L} + I(A_{\cK})$.  
Then there is some 
$$\L \subset \cK ' \subset \cK$$
with $|\cK'| + 1 = |\cK|$ such that $L_{\cK '} \subset Q$.
\end{corollary}
\begin{proof}
We induct on $|\cK| - |\L|$.  If $|\cK| - |\L| > 1$,
we re-index so that $\cK = \{1, \dotsc , s\}$ and 
$\L = \{1, \dotsc , r\}$ with $s > r+1$.  Thus we let
$\J_1 = \{1, \dotsc , r+1\}$ and
$\J_2 = \{1, \dotsc , r,r+2, \dotsc ,s\}.$
Then we can apply Proposition~\ref{LIint}, so either
$L_{\J_1}$ or $L_{\J_2}$ is in $Q$.
If $L_{\J_1}$ is in $Q$ we are done.  If
$L_{\J_2}$ is in $Q$, we are in a smaller case,
and thus done by induction.
\end{proof}

\begin{lemma}\label{Lintlem}
Let $\Delta= \{\J_1,\dotsc,\J_m\}$ be any collection
subsets of $\{1,\dotsc,n\}$ and let $\cK = \cap \J_i$.
Suppose that $Q$ is a prime containing  
$I_{\Delta}(A) + L_{\cK}$.  Then for each $\J_i$ there
is some $\J_j$ such that $Q$ contains 
$L_{\J_i \cap \J_j}$.
\end{lemma}
\begin{proof}
Without loss of generality, let $i=1$.
By Corollary \ref{LIintmp}, there is a $\cK' \supset \cK$
such that $|\cK'| + 1 = |\J_1|$ and $Q$ contains
$L_{\cK'}$.  Re-index so that $\J_1 = \{1\} \cup \cK'$.
Since $\cK' \supset \cK$, there must be
at least one $\J_j$ such that $1 \not\in \J_j$.  Therefore,
$\J_i \cap \J_j \subset \cK'$, so $Q$ contains
$L_{\J_i \cap \J_j}$.
\end{proof}

We now give a lemma which explains the interplay between
$P_{\Delta}$ and $x_{+,\dotsc,+}$.

\begin{lemma}\label{Ppluslem}
Let $\Delta = \{\J,\cK\}$ and $i \in \J \cap \cK$.
Then 
$$(L_{\cK \smallsetminus \{i\}}) \cdot L_{\J} \subset P_{\Delta} + L_{\J \smallsetminus \{i\}}.$$
\end{lemma}
\begin{proof}
Re-index so that $i=1$.  Let $x_{j_1,\dotsc,j_n}$ be any
generator of $L_{\J}$ and $x_{+,k_2,\dotsc,k_n}$ be any
generator of $L_{\cK \smallsetminus \{1\}}$.  Consider the
following element of $J_{\Delta} \subset P_{\Delta}$
$$\det \begin{pmatrix} x_{j_1,\dotsc,j_n}   & x_{j_1,k_2,\dotsc,k_n} \\
                       x_{+,j_2,\dotsc,j_n} & x_{+,k_1,\dotsc,k_n} \end{pmatrix}.$$
Since $x_{+,j_2,\dotsc,j_n} \in L_{\J \smallsetminus \{1\}}$,
the result is clear.
\end{proof}

The next proposition uses the previous results in this section
to show that any minimal prime over $I_{\Delta}$ is made up of
several $P_{\Delta_i}$.  The $\Delta_i$ have the property that
each facet of $\Delta$ is in exactly one $\Delta_i$.

\begin{proposition}\label{equivmp}
Let $\Delta$ be any simplicial complex and $\F_1, \dotsc, \F_m$
its facets.  If $\a$ is any minimal prime containing $I_{\Delta}(A)$,
then there is an equivalence relation on the facets of
$\Delta$, 
$$\F_i \sim \F_j \Longleftrightarrow L_{\F_i \cap \F_j} \not\subset \a.$$
This equivalence relation gives a partition of the facets of
$\Delta$, $\Delta_1 \coprod \dotsb \coprod \Delta_r$ such that
$P_{\Delta_i} \subset \a$ for all $i$.  Moreover, for each $i$,
there is some set $\J \subset \cap \Delta_i$
such that $\J \not\subset \F$ for any facet of $\Delta$ not in
$\Delta_i$, and $\a$ contains $L_{\F \smallsetminus \{j\}}$ for
each $\F \in \Delta_i$ and $j \in \J$.
\end{proposition}
\begin{proof}
It is clear that the relation given is symmetric.  Reflexivity
relies on the minimality of $Q$.  If $Q$ is any prime containing
$I_{\Delta} + L_{\F_i}$ for any facet $\F_i$, then in $T_{\Delta}$,
$Q$ can be expressed as $Q' + L_{\F_i}$ where $Q'$ is an ideal
whose generators are written entirely in terms of the facets
$\F_j$, $j\neq i$.  From this perspecive, it is clear that
$$Q' + \sum_{j\neq i} L_{\F_i \cap \F_j}$$
is also prime, so $Q$ was not minimal.

Transitivity of the relation follows easily from Lemma \ref{Lintlem}.
Suppose that $Q$ contains neither $L_{\F_i \cap \F_j}$ nor
$L_{\F_i \cap \F_k}$.  Then applying Lemma \ref{Lintlem} to the
collection $\{\F_i,\F_j,\F_k\}$, we conclude that $Q$ does not
contain $L_{\F_i \cap \F_j \cap \F_k}$.  Therefore, it cannot
contain $L_{\F_i \cap \F_k}$.

Let $\Delta_i$ be any collection of facets such that for
$\F_j, \F_k \in \Delta_i$, $L_{\F_j \cap \F_k} \not\subset Q.$
Let $\cK = \cap_{\F \in \Delta_i} \F$.  By Lemma \ref{Lintlem},
$L_{\cK} \not\subset Q$.  Therefore, by Propostion~\ref{Qcolon},
$P_{\Delta_i} \subset Q$.

Finally, the last statement is a consequence of the definition
of the equivalence relation, Corollary \ref{LIintmp} and
Lemma \ref{Ppluslem}.
\end{proof}

\subsection{Classification of the Other Minimal Primes}
Next we will show that certain ideals of the kind
mentioned in Proposition~\ref{equivmp} are actually
prime.  If $\Delta$ is a simplicial complex all of
whose facets contain the vertex $k$, let
$\Delta \smallsetminus k$ be the simplicial complex
whose facets are $\J \smallsetminus \{k\}$ for each
facet $\J \in \Delta$.

\begin{thm}\label{classmp}
Let $\Delta$ be any simplicial complex on
$\{1,\dotsc,n\}$ and let 
$\Delta_1 \coprod \dotsb \coprod \Delta_r$
be a partition of the facets of $\Delta$.

For each $\Delta_i$ suppose there is a set
$\cK_i \subset \cap \Delta_i$ such that
for any facet $\J \in \Delta$ which is not 
in $\Delta_i$, $\cK_i \smallsetminus \J$
is nontrivial.  Then
$$\a = \sum_i P_{\Delta_i} + \sum_i \sum_{k \in \cK_i} L_{\Delta_i \smallsetminus k}$$
is a prime ideal.

Any minimal prime over $I_{\Delta}$ has
the form of one of these ideals.
\end{thm}
\begin{proof}
By Proposition~\ref{cov}, we can show this in
$S_{\Delta}$.  Notice that if $\J_i \in \Delta_i$ 
and $\J_j \in \Delta_j$, $i\neq j$, then
$\a \supset L_{\J_i \cap \J_j}$.  Therefore,
$\a$ can be expressed as
$\a_1 + \dotsc + \a_r$ where 
$$\a_i = P_{\Delta_i} + \sum_{k\in\cK_i} L_{\Delta_i,\hat k}$$
is an ideal in $S_{\Delta}$ which is expressed only in terms
of the variables in $S_{\Delta_i}$.  Therefore,
$$S_{\Delta}/\a = S_{\Delta_1}/\a_1 \otimes \dotsb \otimes S_{\Delta_r}/\a_r.$$

Since $\cK_i \subset \Delta_i$, $S_{\Delta_i}/\a_i$ is 
an integral domain for each $i$.  This statement is
true regardless of the field of definition.  Therefore,
$S_{\Delta_i}/\a_i$ remains an integral domain when
it is tensored with the algebraic closure of $\K$.
Thus the tensor product $S_{\Delta}/\a$ is an integral
domain, so $\a$ is prime.

The fact that every minimal prime is of this form
is a consequence of Proposition~\ref{equivmp}.
\end{proof}

\subsection{The Case in Which $\Delta$ Is a Graph}
Now we will give some special cases of Theorem \ref{classmp}.
The first is in the case in which each
facet of $\Delta$ has two elements.  In this case,
$\Delta$ is a graph.

We need one preliminary definition.  For any $j$, let
$$\Delta(j) = \{\J \in \Delta \mid j \in \J\}.$$

\begin{corollary}\label{pairwisemp}
Let $\Delta$ be any graph.  Any minimal prime over
$I_{\Delta}$ is either $P_{\Delta}$ or can be 
expressed as
$$\sum_{j \in \Gamma} P_{\Delta(j)} + \sum_{j\not\in \Gamma} L_{\{j\}}$$
for some vertex cover $\Gamma$ of $\Delta$.
\end{corollary}
\begin{proof}
This is a direct application of Theorem \ref{classmp}.
Since for each facet $\{i,j\}$ of $\Delta$, any prime
containing $I_{\Delta} + \la x_{+,\dotsc,+}\ra$ must
contain either $L_{\{i\}}$ or $L_{\{j\}}$, the statement
about $\Gamma$ being a vertex cover follows.
\end{proof}

\subsection{The Case in Which $\Delta$ Has Two Facets and Our Example}
The second special case we give is when $\Delta$ has
only two facets.

\begin{corollary}\label{2facetmp}
If $\Delta$ is a simplicial complex with
two facets, $\J_1,\J_2$ then the minimal
primes over $I_{\Delta}$ are $P_{\Delta}$ and
$$I_{\Delta} + L_{\J_1 \smallsetminus i_1} + L_{\J_2 \smallsetminus i_2}$$
where $i_1 \not\in \J_2$ and $i_2 \not\in \J_1$.
\end{corollary}

Finally, we give our running example.
\begin{corollary}
Let $\Delta$ have facets $\{\{1,2\},\{1,3\},\{2,3\}\}$.
Then the minimal primes over $I_{\Delta}$ are
\begin{equation*}
\begin{split}
P_{\Delta} & \\
Q_0 &=  I_{\Delta} + L_{\{1\}} + L_{\{2\}} + L_{\{3\}}  \\
Q_1 &=  I_{\{2,3\}} + P_{\{1,2\},\{1,3\}} + L_{\{2\}} + L_{\{3\}} \\
Q_2 &=  I_{\{1,3\}} + P_{\{1,2\},\{2,3\}} + L_{\{1\}} + L_{\{3\}} \\
Q_3 &=  I_{\{1,2\}} + P_{\{1,3\},\{2,3\}} + L_{\{1\}} + L_{\{2\}}
\end{split}
\end{equation*}
unless one of the $a_i=2$, in which case $Q_0$ is not
minimal.
\end{corollary}

\section{Principal Radical Systems and Tables}\label{prssec}

\subsection{Principal Radical Systems In General}
In this section we will show that if $\Delta$ is a
simplicial complex with three or fewer facets, 
$P_{\Delta} = K_{\Delta}+Q_{\Delta}$, which is a prime, 
perfect ideal and $I_{\Delta}(A)$ is radical.
For each of these results we will use principal radical systems.

The notion of a principal radical system has proved very
useful in the study of determinantal ideals.  Hochster and 
Eagon developed it as a method for showing that any ideal 
of minors of a generic matrix was radical.  We follow the 
presentation Bruns and Vetter \cite[\S 12]{brve}.

The main idea is to prove that an ideal is radical by adding
in, one at a time, well-selected elements of the ring until
we have an ideal which is obviously radical.  We will now
cite the theorem as stated in \cite[\S 12]{brve}.

\begin{thm}\label{prs}
Let $R$ be a noetherian ring, and $\F$ a family of ideals in $R$.
Suppose that for every member $I \in \F$ which is not known to
be radical, there is some $x \in R$ such that $I + \la x \ra \in \F$
and one of the following conditions holds:
\begin{enumerate}
\item\label{cond1} $x$ is not a zero-divisor modulo $\rad I$ and $\cap_1^{\infty} (I+ \la x^i \ra)/I = 0$. 
\item\label{cond2} there exists an ideal $J \in \F$, $J \supsetneq I$, such that $xJ \subset I$ and
$x$ is not a zero-divisor modulo $\rad J$.
\end{enumerate}
Then all the ideals $I \in \F$ are radical.
\end{thm}

Note that since all of our rings are graded,
$\cap_1^{\infty} (I+ \la x^i \ra)/I = 0$ will
automatically be satisfied by the Krull Intersection Theorem.
We now apply principal radical systems to the
ideals $P_{\Delta}$, starting with the simplest case, when
$\Delta$ has $1$ facet.

\subsection{The Radicality of $K_{\Delta}+Q_{\Delta}$}
\begin{lemma}\label{1facetrad}
Let $A$ be the generic $a_1 \x \cdots \x a_n$ table and
let $\Gamma$ be any collection of subsets of $\{1,\dotsc,n\}$.
Then $I(A) + L_{\Gamma}$ is radical.
\end{lemma}
\begin{proof}
We induct on $(a_1,\dotsc,a_n)$.  The base case is that in which
$a_i =1$ for all $i$.  In this case, the polynomial ring is
$\K[x_{1,\dotsc,1}]$, which is to say it is has only one variable.
If $L$ is non-empty, then the ideal $L_{\Gamma}$ is generated by
$x_{1,\dotsc,1}$ and if $L$ is empty, the ideal $I(A) + L_{\Gamma}$
is $0$.

For any other $(a_1,\dotsc,a_n)$, consider the following families
of ideals.
\begin{equation*}
\begin{split}
F_{l_1, \dotsc , l_{n}} &= I(A)+ L_{\Gamma} + \la x_{i_1, \dotsc , i_n} \mid (i_1, \dotsc , i_n) \leq_{\text{revlex}}  (l_1, \dotsc , l_n) \ra\\
G_{l_1, \dotsc , l_{n}} &= I(A)+ L_{\Gamma} + \la x_{i_1, \dotsc , i_n} \mid i_j < l_j \text{ for some } j \ra
\end{split}
\end{equation*}

$G_{l_1, \dotsc , l_{n}}$ is radical by induction if any 
$l_i > 1$.  Of course, $G_{1,\dotsc,1} = I(A)+ L_{\Gamma}$.
On the other hand, consider any $l = (l_1, \dotsc , l_{r-1})$.
Let $s(l)$ be the least $l'$ such that $l' > l$.  Let $j$ be 
the least $j$ such that $l_j \neq a_j$.  Then
$$s(l) = (1, \dotsc, 1, l_j + 1, l_{j+1}, \dotsc, l_{r-1}).$$
By definition, $F_l + \la x_{s(l)} \ra = F_{s(l)}$.
Moreover, $G_{s(l)} \supsetneq F_l$ unless 
$l = (a_1, \dotsc, a_{n-1},i)$, in which case $F_l = G_{s(l)}$
and is thus radical.

To show that $x_{s(l)} G_{s(l)} \subset F_l$, let
$x_{i_1, \dotsc , i_n}$ be an arbitrary generator of 
$G_{s(l)}$ which is not contained in $I(A)+ L_{\Gamma}$.  By
the definition of $G_{s(l)}$ there is some $j$ such
that $i_j < s(l)_j$.  By re-indexing, assume $j = 1$ for
ease of notation.  The following minor is in 
$I(A)$
$$\det \begin{pmatrix} x_{s(l)_1, \dotsc, s(l)_n}      & x_{s(l)_1,i_2, \dotsc,i_n} \\
                       x_{i_1,s(l)_2, \dotsc , s(l)_{n}} & x_{i_1, i_2, \dotsc , i_n} \end{pmatrix}.$$
Since $(i_1, s(l)_2, \dotsc ,s(l)_{n}) <_{\text{revlex}} s(l)$,
$(i_1, s(l)_2, \dotsc ,s(l)_{n}) \leq_{\text{revlex}} l$.
Therefore, the antidiagonal product is in $F_l$, and since the
minor is in $I(A) \subset F_l$, the diagonal product is
also in $F_l$.

All that remains to show, then, is that $x_{s(l)}$ is a
nonzerodivisor modulo $\rad G_{s(l)}$.  Since $R/G_{s(l)}$
is isomorphic to $R/(I(A)+L_{\Gamma})$ for smaller values
of the $a_i$, this part is reduced to showing that $x_{1,\dotsc,1}$
is a nonzerodivisor modulo $\rad (I(A) + L_{\Gamma})$.
The minimal primes over $I(A) + L_{\Gamma}$ are $I(A)+L_{\Gamma'}$
where $\Gamma'$ is a collection of subsets of $\{1,\dotsc,n\}$,
each of size $(n-1)$ and such that every set in $\Gamma$ is
contained in a set in $\Gamma'$.  These are prime because
$R/(I(A)+L_{\Gamma'})$ is isomorphic to $R/I(A)$, again for
smaller values of the $a_i$.  Since $x_{1,\dotsc,1}$ is not
in any of the minimal primes, it is a nonzerodivisor modulo
$\rad (I(A) + L_{\Gamma})$.

Therefore, we have shown that $\{F_{l_1, \dotsc , l_{n}},G_{l_1, \dotsc , l_{n}}\}$
is a principal radical system, so $I(A) + L_{\Gamma}$ is
radical.
\end{proof}

This relatively simple case actually is very similar to the more
complicated cases.  We will see very similar arguments again.

\begin{proposition}\label{2facetrad} 
Let $\Delta$ be a simplicial complex with two facets,
$\J_1,\J_2$ and let $\cK$ be a subset of $\J_1 \cup \J_2$.  
Then the ideal $K_{\Delta} + Q_{\Delta} + L_{\cK}$ is radical.
\end{proposition}
\begin{proof}
We re-index so that $\J_1 = \{1,\dotsc,s\}$ and
$\J_2 = \{r,\dotsc,n\}$.

If $\cK$ contains $\J_1$ or $\J_2$, this reduces to
Lemma \ref{1facetrad}, so we suppose that $\cK$ contains
neither $\J_1$ nor $\J_2$.  The minimal primes over
$K_{\Delta} + Q_{\Delta} + L_{\cK}$ are all of the form
$K_{\Delta} + Q_{\Delta} + L_{\Delta, \hat i}$ for some $r \leq i \leq s$
or $K_{\Delta} + Q_{\Delta} + L_{\J_1 \smallsetminus i} +L_{\J_2 \smallsetminus j}$
for some $i < r$ and $j > s$.
This implies that $x_{1,\dotsc,1,+,\dotsc,+}$ is a
nonzero divisor modulo $\rad (K_{\Delta} + Q_{\Delta} + L_{\cK})$.

Consider the following families of ideals.
\begin{equation*}
\begin{split}
F_{l_1, \dotsc , l_{s}} &= K_{\Delta} + Q_{\Delta} + L_{\cK} + \la x_{i_1, \dotsc , i_s,+,\dotsc,+} \mid (i_1, \dotsc , i_s) \leq_{\text{revlex}}  (l_1, \dotsc , l_s) \ra \\
G_{l_1, \dotsc , l_{s}} &= K_{\Delta} + Q_{\Delta} + L_{\cK} + \la x_{i_1, \dotsc , i_n} \mid i_j < l_j \text{ for some } j \leq s \ra
\end{split}
\end{equation*}
The $G_{l_1, \dotsc , l_{s}}$ are defined to allow any of the $i_j = +$,
so long as one of the $i_j$ is a number and $i_j < l_j$.

As in the proof of Lemma \ref{1facetrad}, we can induct on
$(a_1,\dotsc,a_n)$, and thus we can assume that
$G_{l_1, \dotsc , l_{s}}$ as long as one of the $l_i > 1$.
In fact, the entire argument from Lemma \ref{1facetrad}
is valid.  We only need to note that
for any $l$, $x_{s(l)} G_{s(l)} \subset F_l$ and
$F_{a_1,\dotsc,a_s}$ is radical by Lemma \ref{1facetrad}.
\end{proof}

\begin{thm}\label{3facetprime}
Suppose that $\Delta$ is a simplicial complex with no
more than $3$ facets.  Then $K_{\Delta} + Q_{\Delta}$ is radical,
hence $K_{\Delta} + Q_{\Delta} = P_{\Delta}$.
\end{thm}
\begin{proof}
If $\Delta$ has two or fewer facets, then Proposition~\ref{2facetrad}
and Lemma \ref{1facetrad} apply.  Suppose that $\Delta$
has facets $\J_1,\J_2,\J_3$, and re-index so that
$\J_1 = \{1,\dotsc,s\}$.

As in the previous two proofs, consider the following families of
ideals.
\begin{equation*}
\begin{split}
F_{l_1, \dotsc , l_{s}} &= K_{\Delta} + Q_{\Delta} + \la x_{i_1, \dotsc , i_s,+,\dotsc,+} \mid (i_1, \dotsc , i_s) \leq_{\text{revlex}}  (l_1, \dotsc , l_s) \ra \\
G_{l_1, \dotsc , l_{s}} &= K_{\Delta} + Q_{\Delta} + \la x_{i_1, \dotsc , i_n} \mid i_j < l_j \text{ for some } j \leq s \ra
\end{split}
\end{equation*}

They form a principal radical system for the following reasons.  
$G_{l_1, \dotsc , l_{s}}$ is radical by induction on $(a_1,\dotsc,a_n)$.
$F_{1,\dotsc,1}$ satisfies condition~\ref{cond1} of Theorem \ref{prs}
because the radical of $K_{\Delta} + Q_{\Delta}$ is prime.
For $(1,\dotsc,1) \leq l \leq (a_1,\dotsc,a_s)$,
$F_{l}$ satisfies condition~\ref{cond2} of Theorem \ref{prs}
because $F_{s(l)} = F_l + \la x_{s(l),+,\dotsc,+} \ra$ and
$x_{s(l),+,\dotsc,+} \cdot G_{s(l)} \subset F_l$
while $G_{s(l)} \supsetneq F_l$.  Finally,
$F_{a_1,\dotsc,a_s}$ is radical by Proposition~\ref{2facetrad}.

Therefore, we have shown that $K_{\Delta} + Q_{\Delta}$ is 
prime if $\Delta$ has three or fewer faces.
\end{proof}

\subsection{The Perfection of $P_{\Delta}$}
We now use the preceding proofs to establish more about the
algebraic structure of $P_{\Delta}$.  In particular, if $\Delta$
has three or fewer facets, we can show that it is perfect.

\begin{thm}\label{threemain}
If $\Delta$ is a simplicial complex with three or fewer
facets, then $P_{\Delta}$ is a perfect ideal of 
grade $1 - n + \sum a_i$.
\end{thm}
\begin{proof}
We use Proposition~\ref{cov} to reduce to showing that
$P_{\Delta}$ is perfect in the ring $T_{\Delta}$.  Throughout
this proof we will use the same notation as in the previous
proof, and treat all ideals as ideals in $T_{\Delta}$.

The main tool we will use is that if $M_1$, $M_2$, and $M_3$
are $R$-modules such that 
$$0 \lra M_1 \lra M_2 \lra M_3 \lra 0$$
is exact and $M_2$ and $M_3$ are Cohen--Macaulay of depth
$d$ and $d-1$ respectively then
$M_1$ is a Cohen--Macaulay module of depth $d$.

As usual we prove the result by induction on 
$(a_1, \dotsc, a_n)$ since if all but two of these are
$1$, the ideal is just the $2 \x 2$ minors of a generic
matrix, for which this theorem is well-known.

We re-index as in the beginning of the proof of
Theorem \ref{3facetprime}, and use its notation.
By that proof we know that $F_{1, \dotsc , 1}$ 
is radical.  Any prime over $F_{1, \dotsc , 1}$
contains either $P_{\Delta} + L(A_{\J_1})$ or of a
prime of the form
$$H_k = P_{\Delta} + \la x_{i_1, \dotsc , i_n} \mid i_k = 1 \ra$$
for some $k < s$.  Let
$$G_l = \cap_{k = 1}^l H_k.$$

We will show that $R/G_s$ has depth $(\sum a_i)-n$ by
induction.  $R/G_1$ is isomorphic to $R/P_{\Delta}$ with $a_1$
reduced by $1$.  Therefore, $R/G_1 = R/H_1$ is 
Cohen--Macaulay of depth 
$$1 - n + a_1 - 1 + \sum_2^n a_i = (\sum a_i)-n.$$
Now suppose that we have shown that $R/G_k$ has depth
$(\sum a_i)-n$ for any choice of $a_1, \dotsc , a_n$.
Then there is an exact sequence
$$0 \lra R/G_{k+1} \lra R/G_k \oplus R/H_{k+1} \lra R/(G_k + H_{k+1}) \lra 0.$$
the last term is isomorphic to $R/G_k$ where $a_{k+1}$ is
replaced by $a_{k+1}-1$.  Thus, it has depth
$(\sum a_i)-n-1$ by induction.  Both summands of the
middle term have depth $(\sum a_i)-n$ by induction.
Therefore, $R/G_{k+1}$ has $(\sum a_i)-n$.  This implies
that $R/G_s$ is Cohen--Macaulay with depth $(\sum a_i)-n$
as claimed.

If there is some index $j \in \J_1$ but 
$j \not\in \J_l$ for any $l > 1$ then the only minimal primes 
over $F_{1, \dotsc , 1}$ are the $H_k$.  Since $F_{1, \dotsc , 1}$ 
is radical,
we know that $F_{1, \dotsc , 1} = G_r.$  Thus the previous
paragraph implies that $R/F_{1, \dotsc , 1}$ is
Cohen--Macaulay of depth $(\sum a_i)-n$, and since
$$F_{1, \dotsc , 1} = P_{\Delta} + \la x_{1, \dotsc ,1,+,\dotsc,+} \ra$$
and since $P_{\Delta}$ is prime, $x_{1, \dotsc ,1,+,\dotsc,+}$
is a non-zerodivisor modulo it.  Thus $R/P_{\Delta}$
is Cohen--Macaulay of depth $1 - n + \sum a_i$.
Note that if $\Delta$ has two facets (or one), then since neither 
facet can contain the other, this paragraph implies the
theorem for $P_{\Delta}$.

On the other hand, suppose that there is no
$j \in \J_1$ such that $j \not\in \J_l$ for any $l \neq 1$.
This implies that $\Delta$ has three facets,
$\J_1,\J_2,\J_3$.  We may assume that the condition
holds for $\J_2,\J_3$ as well, so for each $j \in \{1,\dotsc,n\}$,
$j$ is an element of two of $\J_1,\J_2,\J_3$.
Therefore, $\J_1$ must contain the symmetric difference
of $\J_2$ and $\J_3$, $(\J_2 \cup \J_3) \smallsetminus (\J_2 \cap \J_3)$.
Thus the minimal primes over $P_{\Delta} + L_{\J_1}$ are 
$$D_i = P_{\Delta} + L_{\J_1} + L_{\Delta,\hat i}$$
for each $i$ in $(\J_2 \cap \J_3) \smallsetminus \J_1$.  The $D_i$
are prime because if $i \in (\J_2 \cap \J_3) \smallsetminus \J_1$,
$R/D_i$ is isomorphic to $R/(P_{\J_2,\J_3})$ with $a_i$ reduced by $1$.
Thus, these prime ideals are also perfect of grade $(\sum a_i)-n$ by
induction.  Our next goal is to prove that their intersection
is also perfect.

Re-index so that
$$(\J_2 \cap \J_3) \smallsetminus \J_1 = \{1,\dotsc,m\}$$
and let 
$$E_l = \cap_{i=1}^l D_i.$$
Suppose that $E_k$ is perfect of grade $(\sum a_i)-n$.
Then we have an exact sequence
$$0 \lra R/E_{k+1} \lra R/E_k \oplus R/D_{k+1} \lra R/(E_k + D_{k+1}) \lra 0.$$
We know that $R/E_k$ and $R/D_{k+1}$ are both
Cohen--Macaulay of depth $(\sum a_i)-n$, and
since $R/(E_k + D_{k+1}) \cong (R/D_{k+1})/E_k$ which
is isomorphic to $R/E_k$ with $a_{k+1}$ decreased by
$1$, $R/(E_k + D_{k+1})$ is Cohen--Macaulay of depth
$(\sum a_i)-n-1$.  Therefore, we know that $R/E_{k+1}$ 
is Cohen--Macaulay of depth $(\sum a_i)-n$.
Therefore, by induction, $(P_{\Delta} + L_{\J_1})$
is perfect of grade $(\sum a_i)-n$.

Finally, we need to show that
$$F_{1, \dotsc , 1} = (P_{\Delta} + L_{\J_1}) \cap G_s,$$
is also perfect of grade $(\sum a_i)-n$, where 
$G_s = \cap H_k$ is defined as above.  This can be established
in exactly the same way as the perfection of $G_s$ and
$P_{\Delta} + L_{\J_1}$ were.  Let $C_k = (P_{\Delta} + L_{\J_1}) \cap G_k$,
where $G_0 = \la 1 \ra$.
We have already established that $P_{\Delta} + L_{\J_1}$
is perfect of grade $(\sum a_i)-n$, so suppose that
$C_k$ is perfect.  We have $C_{k+1} = C_k \cap H_{k+1}$ and
thus an exact sequence
$$0 \lra R/C_{k+1} \lra R/C_k \oplus R/H_{k+1} \lra R/(C_k + H_{k+1}) \lra 0.$$
Like the previous proofs, $R/C_k$ and $R/H_{k+1}$ we already
know to be Cohen--Macaulay of depth $(\sum a_i)-n$, and
$R/(C_k + H_{k+1}) \cong (R/H_{k+1})/C_k$, which is isomorphic
to $R/C_k$ for $a_{k+1}$ decreased by $1$, so it is Cohen--Macaulay
of depth $(\sum a_i)-n-1$.  Therefore, $R/C_{k+1}$ is Cohen--Macaulay
of depth $(\sum a_i)-n$, so by induction,
$R/C_s = R/F_{1, \dotsc , 1}$ is Cohen--Macaulay of depth $(\sum a_i)-n$.

Since $F_{1, \dotsc , 1} = P_{\Delta} + \la x_{1,\dotsc,1,+,\dotsc,+} \ra$
and $x_{1,\dotsc,1,+,\dotsc,+}$ is a nonzerodivisor modulo $P_{\Delta}$,
this implies that $P_{\Delta}$ is perfect of grade $1-n + \sum a_i$.
\end{proof}

\subsection{The Radicality of $I_{\Delta}$}
We now move from the prime ideal $P_{\Delta}$
to the original ideal $I_{\Delta}$.

\begin{proposition}\label{2facetIrad}
Let $\Delta$ be a simplicial complex with two facets,
$\J_1,\J_2$ and let $\cK$ be a subset of $\J_1 \cup \J_2$.  
Then the ideal $I_{\Delta} + L_{\cK}$ is radical.
\end{proposition}
\begin{proof}
We re-index so that $\J_1 = \{1,\dotsc,s\}$ and
$\J_2 = \{r,\dotsc,n\}$.

If $\cK$ contains $\J_1$ or $\J_2$, this reduces to
Lemma \ref{1facetrad}, so we suppose that $\cK$ contains
neither $\J_1$ nor $\J_2$.  We will prove the result
by pricipal radical systems.  Define 
\begin{equation*}
\begin{split}
F_{l_1, \dotsc , l_s} &= I_{\Delta} + L_{\cK} + \la x_{i_1, \dotsc , i_s,+,\dotsc,+} \mid (i_1, \dotsc , i_s) \leq_{\text{revlex}}  (l_1, \dotsc , l_s) \ra \\
G_{l_1, \dotsc , l_s} &= I_{\Delta} + L_{\cK} + \la x_{i_1, \dotsc , i_s,+,\dotsc,+} \mid i_j < l_j \text{ for some } j \leq s \ra
\end{split}
\end{equation*}
and let $\F = \{F_{l_1, \dotsc , l_s} + G_{k_1, \dotsc , k_s}\}$,
the set of all sums of $F$'s and $G$'s.  We claim that
$\F$ is a principal radical system.

If $l = (l_1, \dotsc , l_s)$ is any sequence, let $s(l)$ 
be the least $l'$ such that $l' >_{\text{revlex}} l$.  
If $j$ is the least $j$ such that $l_j \neq a_j$ then
$$s(l) = (1, \dotsc, 1, l_j + 1, l_{j+1}, \dotsc, l_s).$$
By definition, $F_l + \la x_{s(l)} \ra = F_{s(l)}$.
Therefore, $F_l + G_k + \la x_{s(l)} \ra = F_{s(l)} + G_k$.

The following lemma will be the key to showing that
$\F$ is a principal radical system.
\begin{lemma}\label{2facetlem}
$x_{l_1,\dotsc,l_s,+,\dotsc,+}$ is a nonzerodivisor
modulo $\rad G_{l_1,\dotsc,l_s}$.
\end{lemma}
\begin{proof}
We will do this by computing the minimal primes
over $G_{l_1,\dotsc,l_s}$, and showing that
$x_{l_1,\dotsc,l_s,+,\dotsc,+}$ is not in any of them.
Let $l' = (1,\dotsc,1,l_r,\dotsc,l_s)$.  Then
$R/G_l \cong R/G_{l'}$ where the latter ring has the
values of $a_i$ decreased by $l_i-1$ for each $i < r$.
Therefore, we can assume that $l_i =1$ for all $i<r$.

Suppose that $l_i > 1$ for some $i \geq r$, without
loss of generality, assume $i = s$.  Then for each
$j_r,\dotsc,j_{s-1}$ and any $j_s < l_s$
$$x_{+,\dotsc,+,j_r,\dotsc,j_s,+,\dotsc,+} \in G_l.$$
Since $I(A_{\J_2}) \subset G_l$, any prime
containing $G_l$ must either contain $L_{\{r,\dotsc,s\}}$
or $x_{+,\dotsc,+,j_r,\dotsc,j_s,j_{s+1},\dotsc,j_n}$
for all $j_s < l_s$.
 
Let $H_{l_1,\dotsc,l_s} = G_{l_1,\dotsc,l_s} + \la x_{+,\dotsc,+,i_r, \dotsc , i_n} \mid i_j < l_j \text{ for some } r \leq j \leq s \ra.$
The previous paragraph implies that any prime
containing $G_{l_1,\dotsc,l_s}$ either contains
$L_{\{r,\dotsc,s\}}$ or contains $H_{l_1,\dotsc,l_s}$.
Since $R/H_l$ is isomorphic to $R/(I_{\Delta}+L_{\cK})$
with $a_i$ decreased by $l_i-1$ for each $i$.  Therefore,
to show that $x_{l_1,\dotsc,l_s,+,\dotsc,+}$ is not
in a minimal prime over $H_{l_1,\dotsc,l_s}$ is the same
as showing that $x_{1,\dotsc,1,+,\dotsc,+}$ is not in
a minimal prime over $I_{\Delta} + L_{\cK}$.

The minimal primes over $I_{\Delta} + L_{\cK}$ are
either $P_{\Delta} + L_{\Delta,\hat i}$ for some
$i \in (\J_1 \cap \J_2) \smallsetminus \cK$ or
$I_{\Delta} + L_{\J_1 \smallsetminus i_1} + L_{\J_2 \smallsetminus i_2}$
where $i_1,i_2 \not\in \cK$.  It is clear that
$x_{1,\dotsc,1,+,\dotsc,+}$ is not in any of these
ideals.

On the other hand, we must show that 
$x_{1,\dotsc,1,+,\dotsc,+}$ is not in any of
the minimal primes over $G_l + L_{\{r,\dotsc,s\}}$.
Because this ideal contains $L_{\{r,\dotsc,s\}}$,
it can be expressed, in $S_{\Delta}$ as $I_1 + I_2$
where $I_1 \subset \K[X_{i_1,\dotsc,i_s,+,\dotsc,+}]$ 
and $I_2 \subset \K[X_{+,\dotsc,+,i_r,\dotsc,i_n}]$.
Therefore, we need only consider the minimal primes over
$$I(A_{\J_1}) + L_{\cK\cap \J_1} + L_{\{r,\dotsc,s\}} + \la x_{i_1, \dotsc , i_s,+,\dotsc,+} \mid i_j < l_j \text{ for some } j \leq s \ra.$$
The effect of the last summand is only to reduce
each $a_i$ by $l_i-1$, so we may assume that this
term is $0$.  Then we are left with 
$I(A_{\J_1}) + L_{\cK\cap \J_1} + L_{\{r,\dotsc,s\}}$,
whose minimal primes are contained in
$I(A_{\J_1}) + L_{\J_1\smallsetminus i} + L_{\J_1\smallsetminus j}$
where $i \not\in \cK$ and $j <r$.  Thus
$x_{1,\dotsc,1,+,\dotsc,+}$ is not in any minimal
prime over $G_l+L_{\{r,\dotsc,s\}}$.

This completes the proof of the lemma, so
$x_{l,+,\dotsc,+}$ is a nonzerodivisor modulo
$G_l$.
\end{proof}

Since $G_{s(l)} \supsetneq F_l$ for any $l$,
$F_l + G_k = G_k$ whenever $k >_{\text{revlex}} l$.
Moreover,
$$G_{l} + F_{l} = G_{l} + \la x_{l_1,\dotsc,l_s,+,\dotsc,+} \ra$$
so by our lemma, if $k >_{\text{revlex}} l$ 
$F_l + G_k$ satisfies condition~\ref{cond1} of
theorem \ref{prs}.

On the other hand, if $k \leq l < (a_1,\dotsc,a_n)$, recall that
$$G_k + F_{s(l)} = G_k + F_l + \ra x_{s(l),+,\dotsc,+}.$$
$G_{s(l)} \supsetneq G_k + F_l$, and 
$x_{s(l),+,\dotsc,+} G_{s(l)} \subset G_k + F_l$.
Thus, since $x_{s(l),+,\dotsc,+}$ is a nonzerodivisor
modulo $\rad G_{s(l)}$ by the lemma, $G_k + F_l$
satisfies condition~\ref{cond2} of theorem \ref{prs}.

Finally, 
$$F_{a_1,\dotsc,a_n}= I_{\Delta} + L_{\cK} + L_{\J_1} = I(A_{\J_2}) + L_{\cK} + L_{\J_1}$$
which is radical by Lemma \ref{1facetrad}.

Therefore, $\F$ is a principal radical system and
$I_{\Delta} + L_{\cK}$ is radical.
\end{proof}

\begin{thm}
If $\Delta$ has three or fewer facets then
$I_{\Delta}$ is a radical ideal.
\end{thm}
\begin{proof}
If $\Delta$ has one or two facets, this has been
proven in Lemma \ref{1facetrad} and Proposition
\ref{2facetIrad}, so we may assume that $\Delta$
has three facets.

This proof is very similar to that of Proposition
\ref{2facetIrad}.  Re-index so that $\J_1 = \{1,\dotsc,s\}$,
and let
\begin{equation*}
\begin{split}
F_{l_1, \dotsc , l_s} &= I_{\Delta} + \la x_{i_1, \dotsc , i_s,+,\dotsc,+} \mid (i_1, \dotsc , i_s) \leq_{\text{revlex}}  (l_1, \dotsc , l_s) \ra \\
G_{l_1, \dotsc , l_s} &= I_{\Delta} + L_{\cK} + \la x_{i_1, \dotsc , i_s,+,\dotsc,+} \mid i_j < l_j \text{ for some } j \leq s \ra.
\end{split}
\end{equation*}
Define $\F = \{F_{l_1, \dotsc , l_s} + G_{k_1, \dotsc , k_s}\}$,
the set of all sums of $F$'s and $G$'s.  We claim that
$\F$ is a principal radical system.

We will prove below that $x_{l,+,\dotsc,+}$
is a nonzerodivisor modulo $\rad G_l$ and
now we show how that will imply the theorem.

As in the previous proof, 
$F_l + G_k + \la x_{s(l)} \ra = F_{s(l)} + G_k$,
so so for if $k >_{\text{revlex}} l$ 
$F_l + G_k$ satisfies condition~\ref{cond1} of
theorem \ref{prs}.  Moreover, if 
$k \leq l < (a_1,\dotsc,a_n)$, 
$$G_k + F_{s(l)} = G_k + F_l + \ra x_{s(l),+,\dotsc,+}.$$
$G_{s(l)} \supsetneq G_k + F_l$, and 
$x_{s(l),+,\dotsc,+} G_{s(l)} \subset G_k + F_l$.
Thus, since $x_{s(l),+,\dotsc,+}$ is a nonzerodivisor
modulo $\rad G_{s(l)}$ as we will show below, $G_k + F_l$
satisfies condition~\ref{cond2} of theorem \ref{prs}.

Finally, $F_{a_1,\dotsc,a_s}$ is radical because it is
$I_{\J_2,\J_3} + L_{\J_1}$ which is radical by
Proposition~\ref{2facetIrad}.

Therefore, the theorem will be completed with the proof
of the following lemma
\begin{lemma}
$x_{l_1,\dotsc,l_s,+,\dotsc,+}$ is a nonzerodivisor modulo
$\rad G_{l_1,\dotsc,l_s}$.
\end{lemma}
\begin{proof}
Again, we prove this by computing the minimal primes over
$G_{l_1,\dotsc,l_s}$ and showing that $x_{l_1,\dotsc,l_s,+,\dotsc,+}$
is not in any of them.

Suppose $Q$ is a minimal prime over $G_{l_1,\dotsc,l_s}$ which
contains $L_{\J_1\cap\J_2} + L_{\J_1\cap\J_3}$. Since
$G_{l_1,\dotsc,l_s} + L_{\J_1\cap\J_2} + L_{\J_1\cap\J_2}$ can
be expressed as $I_1 + I_2$ where 
$I_1 \subset \K[X_{i_1,\dotsc,i_s,+,\dotsc,+}] = S_{\J_1}$ and
the generators $I_2$ have none of those variables in them,
we can show that $x_{l_1,\dotsc,l_s,+,\dotsc,+} \not\in Q$
by showing that it is not in any minimal prime over
$I(A_{\J_1}) + L_{\J_1\cap\J_2} + L_{\J_1\cap\J_3}$, which
is clear.

The second case is when $Q$ is a minimal prime over 
$G_{l_1,\dotsc,l_s}$ which contains $L_{\J_1\cap\J_3}$
but not $L_{\J_1\cap\J_2}$.  Then it must also contain
$L_{\J_2 \cap \J_3}$ and $P_{\{\J_1,\J_2\}}$ by Corollary 
\ref{equivmp}.  As in the previous paragraph, 
$$G_{l_1,\dotsc,l_s} + P_{\J_1,\J_2} + L_{\J_1\cap\J_3} + L_{\J_2\cap\J_3}$$
can be expressed as $I_1 + I_2$ where 
$I_1 \subset S_{\{\J_1,\J_2}\}$ and the generators $I_2$ 
have none of those variables in them.  Thus we can show that
$x_{l_1,\dotsc,l_s,+,\dotsc,+} \not\in Q$ by showing that
it is not in any minimal prime over
$P_{\{\J_1,\J_2\}} + G_{l_1,\dotsc,l_s} + L_{\J_3}.$
Since this case was covered in Lemma \ref{2facetlem},
we refer to that proof.

The final case is that in which $Q$ is a minimal prime over 
$G_{l_1,\dotsc,l_s}$ and contains neither $L_{\J_1\cap\J_2}$
nor $L_{\J_1\cap\J_3}$.  Thus, it cannot contain $L_{\J_1\cap\J_3}$
either and must contain $P_{\Delta}$.  Moreover, if 
$i \in \J_1 \cap \J_2$ and $l_i > 1$, then we can re-index so
$i=s$ and $\J_2 = \{r,\dotsc,n\}$.  As in the proof of Lemma
\ref{2facetlem}, $G_{l_1,\dotsc,l_s}$ contains
$x_{+,\dotsc,+,j_r,\dotsc,j_s,+,\dotsc,+}$
for all $j_r,\dotsc,j_{s-1}$ and any $j_s < l_s$.
Since $Q$ does not contain $L_{\J_1 \cap \J_2}$
it must be the case that
$x_{+,\dotsc,j_r,\dotsc,j_n} \in Q$ as long
as $j_s < l_s$.  Therefore, $Q$ must contain
$$H_{l_1,\dotsc,l_s} = P_{\Delta} + \la x_{i_1, \dotsc , i_n} \mid i_j < l_j \text{ for some } j \leq s \ra.$$
(Notice that this ideal was defined as $G_{l_1,\dotsc,l_s}$ in the
proof of Theorem \ref{3facetprime}.)  Since
$R/H_{l_1,\dotsc,l_s}$ is isomorphic to
$R/P_{\Delta}$ where each $a_i$ has been reduced by
$l_i - 1$.  Therefore, $H_{l_1,\dotsc,l_s}$ is
prime and $x_{l_1,\dotsc,l_s,+,\dotsc,+}$ is not
in it.

We have shown that $x_{l_1,\dotsc,l_s,+,\dotsc,+}$
is not in any minimal prime over
$G_{l_1,\dotsc,l_s}$ and hence is a nonzerodivisor
modulo its radical.
\end{proof}
This completes the proof that $I_{\Delta}$ is radical if $\Delta$
has fewer than three facets.
\end{proof}

\section{Conjectures, Examples, and Notes on Computation}\label{finalsec}
\subsection{An Example In Which $I_{\Delta}$ Is Not Radical}

It is not true that for any $\Delta$, $I_{\Delta}(A)$
is radical.  Any time $Q_{\Delta} \neq J_{\Delta}$, we know that
$$x_{+,\dotsc,+} \cdot Q_{\Delta} \subset \rad(I_{\Delta});$$
however, this will not always be contained in $I_{\Delta}$.
For example, when 
$$\Delta = \{\{1,2\},\{1,3\},\{2,4\},\{3,4\}\}$$
it can be shown computationally that
$$x_{+,+,+,+}(x_{1,1,+,+}x_{+,+,1,1}-x_{1,+,1,+}x_{+,1,+,1}) \not\in I_{\Delta}.$$

In this case, it turns out that the primary decomposition
is still accessible, and we give a computation of it in the
case $a_i=2$.  Let 
\begin{equation*}
\begin{split}
Q_1 &= P_{\{1,2\},\{2,4\}} + P_{\{1,3\},\{3,4\}} + \la x_{i,+,+,+} \ra + \la x_{+,+,+,l} \ra \\
Q_2 &= P_{\{1,2\},\{1,3\}} + P_{\{2,4\},\{3,4\}} + \la x_{+,j,+,+} \ra + \la x_{+,+,k,+} \ra \\
Q_3 &= I_{\Delta} + \la x_{i,+,+,+}^2, x_{+,j,+,+}^2, x_{+,+,k,+}^2, x_{+,+,+,l}^2,x_{+,+,+,+}^2 \ra \\
\end{split}
\end{equation*}
It can be verified using Macaulay 2 \cite{M2} that
$$I_{\Delta} = P_{\Delta} \cap Q_1 \cap Q_2 \cap Q_3.$$

\subsection{Two Conjectures}
Section~\ref{prssec} has exclusively dealt with the case 
in which $\Delta$ has three or fewer facets.  We offer the following
conjectures which have been borne out in all the examples
which our computers have been able to accomplish.

\begin{conjecture}
If $\Delta$ is any simplicial complex,
$$K_{\Delta} + Q_{\Delta} = P_{\Delta},$$ 
which is a prime and perfect ideal of grade
$1-n+\sum a_i$.
\end{conjecture}

We have proven this result in the case in which $\Delta$
has three or fewer facets.  Moreover, we have shown that
$\rad (K_{\Delta}+Q_{\Delta})$ is prime in Theorem \ref{radsegeq}, 
which should be seen as good evidence for the primality 
of the ideal.

The second conjecture deals with the radicality of $I_{\Delta}$.
\begin{conjecture}
Let $\Delta$ be any simplicial complex.  $I_{\Delta}$
is a radical ideal if and only if $Q_{\Delta} = J_{\Delta}$.
\end{conjecture}

\subsection{Notes on Computation}
Finally, we discuss the the copmutational aspects of experimenting
with these families of ideals.  All computations
should be done in $T_{\Delta}$ because it reduces the number of
variables in the polynomial ring.  This reduction is especially
noticeable when some of the $a_i > 2$.  A side benefit is that
the relations are usually easier to decipher when they are
expressed in the variables of $S_{\Delta}$.  In fact, these were
the reasons that first attracted me to change variables.

If $a_i =2$ for all $i$, then we are in a position to decompose 
$I_{\Delta}$ when $\Delta$ has fewer than four vertices 
($n \leq 4$), and can do some cases with five or six vertices.
After that point, the only $\Delta$'s for which $I_{\Delta}$ can 
be decomposed have two facets.

When $a_i = 2$ it is also possible to compute a free resolution
for $P_{\Delta}$ for some cases until $n = 5$.  After that,
the problem again becomes insurmountable.

If we allow $a_i >2$, both problems become very difficult very fast.
The decomposition can be checked by using Theorem \ref{classmp},
and intersecting the minimal primes.  Computing a free resolution
also becomes computationally impossible very fast.  For
the simplest $\Delta$ with three facets, $\{\{1,2\},\{1,3\},\{2,3\}\}$,
a free resolution cannot be computed when $a_i =3$ for each $i$.

\providecommand{\bysame}{\leavevmode\hbox to3em{\hrulefill}\thinspace}
\providecommand{\xxMR}[1]{\relax\ifhmode\unskip\space\fi
  \href{http://www.ams.org/mathscinet-getitem?mr=#1}{MR #1}}
\providecommand{\xxZBL}[1]{\relax\ifhmode\unskip\space\fi
  \href{http://www.emis.de/cgi-bin/zmen/ZMATH/en/quick.html?type=html&an=#1}{Z%
BL=#1}}
\providecommand{\xxARXIV}[1]{\relax\ifhmode\unskip\space\fi
  \href{http://arxiv.org/abs/#1}{arXiv=#1}}
\providecommand{\href}[2]{#2}


\begin{thebibliography}{Mat99}

\bibitem[BV88]{brve}
Winfried Bruns and Udo Vetter, \emph{Determinantal rings}, Lecture Notes in
  Mathematics, vol. 1327, Springer-Verlag, Berlin, 1988. \xxMR{89i:13001}

\bibitem[GS]{M2}
Daniel~R. Grayson and Michael~E. Stillman,
  \href{http://www.math.uiuc.edu/Macaulay2/}{\emph{Macaulay 2, a software
  system for research in algebraic geometry}}.

\bibitem[H{\`a}02]{HTHa}
Huy~T{\`a}i H{\`a},
  \href{http://dx.doi.org/10.1016/S0022-4049(01)00032-9}{\emph{Box-shaped
  matrices and the defining ideal of certain blowup surfaces}}, J. Pure Appl.
  Algebra \textbf{167} (2002), no.~2-3, 203--224. \xxMR{2002h:13020}

\bibitem[Mat99]{Matus}
F.~Mat{\'u}{\v{s}},
  \href{http://journals.cambridge.org/bin/bladerunner?REQUNIQ=1095195553&REQSE%
SS=3451300&118200REQEVENT=&REQINT1=46684&REQAUTH=0}{\emph{Conditional
  independences among four random variables. {III}. {F}inal conclusion}},
  Combin. Probab. Comput. \textbf{8} (1999), no.~3, 269--276.
  \xxMR{2000i:68176}

\bibitem[Stu02]{SPE}
Bernd Sturmfels, \emph{Solving systems of polynomial equations}, CBMS Regional
  Conference Series in Mathematics, vol.~97, Published for the Conference Board
  of the Mathematical Sciences, Washington, DC, 2002. \xxMR{2003i:13037}

\end{thebibliography}
\end{document}